\documentclass[aop,MSNbibl,seceqn,dvips]{arximspdf}

\doi{10.1214/12-AOP760} %kopijuoti is PTS
\volume{41}
\issue{5}
\pubyear{2013}
\firstpage{3081}
\lastpage{3111}

\makeatletter

\newcommand{\rref}[1]{\mbox{\fontsize{8.36pt}{9pt}\selectfont{\ref{#1}}}}

\newtheorem{theorem}{Theorem}[section]
\newtheorem{proposition}[theorem]{Proposition}
\newtheorem{corollary}[theorem]{Corollary}
\newtheorem{lemma}[theorem]{Lemma}

\newproclaim{remark}{Remark}

\newcommand{\E}{\mathbb{E}}
\newcommand{\range}{\operatorname{range}}
\newcommand{\rank}{\operatorname{rank}}
\newcommand{\Var}{\operatorname{Var}}
\newcommand{\tr}{\operatorname{tr}}

\newcommand{\R}{\mathbb{R}}
\renewcommand{\P}{\mathbb{P}}
\newcommand{\one}{{\mathbf1}}
\renewcommand{\l}{\lambda}
\renewcommand{\t}{\tau}

\newcommand{\dn}{\delta}
\newcommand{\D}{\Delta}
\newcommand{\Qn}{q}
\newcommand{\Qp}{Q}
\newcommand{\mn}{\underline{m}}
\renewcommand{\mp}{\overline{m}}
\newcommand{\pn}{\phi}
\newcommand{\lmin}{\lambda_{\mathrm{min}}}
\newcommand{\lmax}{\lambda_{\mathrm{max}}}
\newcommand{\SR}{\texttt{SR}}
\newcommand{\upper}{{\mathrm{upper}}}
\newcommand{\low}{{\mathrm{lower}}}
\newcommand{\main}{{\mathrm{main}}}
\newcommand{\prop}{{\mathrm{trace}}}

\makeatother

\begin{document}
\begin{frontmatter}

\title{Covariance estimation for distributions with~$2+\varepsilon$~moments}
\runtitle{Covariance estimation}

\begin{aug}
\author[A]{\fnms{Nikhil} \snm{Srivastava}\corref{}\thanksref{t1}\ead[label=e1]{nikhils@math.ias.edu}}
\and
\author[B]{\fnms{Roman} \snm{Vershynin}\thanksref{t2}\ead[label=e2]{romanv@umich.edu}}
\runauthor{N. Srivastava and R. Vershynin}
\affiliation{Institute for Advanced Study and University of Michigan}
\address[A]{School of Mathematics\\
Institute for Advanced Study\\
1 Einstein Drive\\
Princeton, New Jersey 08540\\
USA\\
\printead{e1}} %adresu isvedimo komanda gale!
\address[B]{Department of Mathematics\\
University of Michigan\\
530 Church St.\\
Ann Arbor, Michigan 48109\\
USA\\
\printead{e2}}
\end{aug}

\thankstext{t1}{Supported by NSF Grants CCF 0832797 and DMS-08-35373.}

\thankstext{t2}{Supported in part by NSF Grants FRG DMS-09-18623 and
DMS-10-01829.}

% HISTORY:
\received{\smonth{6} \syear{2011}}
\revised{\smonth{3} \syear{2012}}

% ABSTRACT
%
\begin{abstract}
We study the minimal sample size $N=N(n)$ that suffices to estimate the
covariance matrix of an $n$-dimensional distribution by the sample
covariance matrix in the operator norm, with an arbitrary fixed
accuracy. We establish the optimal bound $N = O(n)$ for every
distribution whose $k$-dimensional marginals have uniformly bounded
$2+\varepsilon$ moments outside the sphere of radius $O(\sqrt{k})$. In
the specific case of log-concave distributions, this result provides an
alternative approach to the Kannan--Lovasz--Simonovits problem, which
was recently solved by Adamczak et al. [\textit{J. Amer. Math. Soc.}
\textbf{23} (2010) 535--561]. Moreover, a lower estimate on the
covariance matrix holds under a weaker assumption---uniformly bounded
$2+\varepsilon$ moments of one-dimensional marginals. Our argument
consists of randomizing the spectral sparsifier, a~deterministic tool
developed recently by Batson, Spielman and Srivastava [\textit{SIAM J.
Comput.} \textbf{41} (2012) 1704--1721]. The new randomized method
allows one to control the spectral edges of the sample covariance
matrix via the Stieltjes transform evaluated at carefully chosen random
points.
\end{abstract}

% KEYWORDS
%
\begin{keyword}[class=AMS]
\kwd[Primary ]{62H12}
\kwd[; secondary ]{60B20}
\end{keyword}
\begin{keyword}
\kwd{Covariance matrices}
\kwd{high-dimensional distributions}
\kwd{Stieltjes transform}
\kwd{log-concave distributions}
\kwd{random matrices}
\end{keyword}

\end{frontmatter}

%s1 #&#
\section{Introduction}\label{intro}

%s1.1 #&#
\subsection{Covariance estimation problem}

Estimating covariance matrices of high-dimensional distributions is a
basic problem
in statistics and its numerous applications. Consider a random vector
$X$ valued in $\R^n$,
and let us assume for simplicity that $X$ is centered, that is, $\E X = 0$;
this restriction will not be needed later.
The covariance matrix of $X$ is the $n \times n$ positive semidefinite matrix
\[
\Sigma= \E X X^T.
\]
Our goal is to estimate $\Sigma$ from a sample $X_1, \ldots, X_N$
taken from the
same distribution as $X$. A classical unbiased estimator for $\Sigma$ is
the sample covariance matrix
\[
\Sigma_N = \frac{1}{N} \sum_{i=1}^N
X_i X_i^T.
\]
A basic question is to \textit{determine the minimal sample size $N$ which
guarantees that $\Sigma$ is accurately estimated by $\Sigma_N$.}
More precisely, for a given accuracy $\varepsilon> 0$, we are interested
in the minimal $N = N(n,\varepsilon)$ so that
\[
\E\|\Sigma_N - \Sigma\| \le\varepsilon\|\Sigma\|,
\]
where \mbox{$\|\cdot\|$}\vspace*{1pt} denotes the spectral (operator) norm.
Replacing $X$ by $\Sigma^{-1/2} X$ and $X_i$ by $\Sigma^{-1/2} X_i$,
we reduce the problem to the distributions for which $\Sigma= I$, that
is, to \textit{isotropic distributions}.

%s1.2 #&#
\subsection{Sampling from isotropic distributions}

We consider independent iso\-tropic random vectors $X_i$ valued in $\R^n$, that is, such that
$\E X_i X_i^T = I$.
Our goal is to determine the minimal sample size $N = N(n,\varepsilon
)$ such that
\[
\E\|\Sigma_N - \Sigma\| \le\varepsilon.
\]
For obvious-dimensional reasons, one must have $N \ge n$.
Rudelson's remarkably general result (\cite{R}, see~\cite{Vtutorial},
Section 4.3) yields that
if $\|X\|_2 = O(\sqrt{n})$ almost surely, then
%
%e1.1 #&#
\begin{equation}
\label{nlogn} N = O(n \log n),
\end{equation}
where the $O(\cdot)$ notation hides the dependence on $\varepsilon$
here and thereafter.
It is well known that the logarithmic oversampling factor cannot be
removed from
(\ref{nlogn}) in general, %can not->cannot
for example, if the distribution is supported on $O(n)$ points; see
Section~\ref{soptimality}.

Nevertheless, it is also known that for sufficiently regular
distributions the logarithmic oversampling factor is not
needed in (\ref{nlogn}). This is a property of the standard normal
distribution in $\R^n$ and, more generally, of the distributions with
\textit{sub-Gaussian}
one-dimensional marginals. Namely,
\[
N = O(n)
\]
holds for every distribution that satisfies
%
%e1.2 #&#
\begin{equation}
\label{sub-Gaussian} \sup_{\|x\|_2 \le1} \bigl(\E\bigl|\langle X, x
\rangle\bigr|^p\bigr)^{1/p} = O\bigl(\sqrt {p}\bigr) \qquad\mbox{for $p \ge1$}.
\end{equation}
This result can be obtained by a standard covering argument; see
\cite{Vtutorial}, Section 4.3.

It is an open problem to describe the distributions for which the
logarithmic oversampling
is not needed, that is, for which $N = O(n)$.
The gap between sub-Gaussian distributions where this bound holds and
discrete distributions
on $O(n)$ points where it fails is quite large. %really big.

It is already a difficult problem to relax the sub-Gaussian moment
assumption (\ref{sub-Gaussian})
to anything weaker while keeping $N = O(n)$.
A major step was made by
Adamczak et al.~\cite{ALPT}, who
showed that $N = O(n)$ still holds (in fact, with high probability)
under the \textit{sub-exponential} moment assumptions
%
%e1.3 #&#
\begin{eqnarray}
\label{sub-exponential} \|X\|_2 &=& O\bigl(\sqrt{n}\bigr)
\qquad\mbox{a.s.},\nonumber\\[-8pt]\\[-8pt]
\sup_{\|x\|_2 \le1} \bigl(\E\bigl|\langle X, x\rangle\bigr|^p
\bigr)^{1/p} &=& O(p) \qquad\mbox{for $p \ge1$.}\nonumber
\end{eqnarray}
As an application, it was shown in~\cite{ALPT} that $N = O(n)$ holds
for log-concave distributions,
and in particular for the uniform distributions on isotropic convex
bodies in $\R^n$. This answered a
question posed by Kannan, Lovasz and Simonovits in~\cite{KLS}.

The second author of the present paper speculated in~\cite{Vcovariance}
that $N = O(n)$ should hold for a much wider class of distributions
than sub-exponential, perhaps \textit{for all distributions with
$2+\varepsilon$ moments}. (The second moment---the variance---is
assumed to be finite by the nature of the problem, as otherwise the
covariance matrix is not defined.) The goal of the the current paper is
to provide a result of this type.

%th1.1 #&#
\begin{theorem} \label{main}
Consider independent isotropic random vectors $X_i$ valued in~$\R^n$.
Assume that $X_i$ satisfy the \textup{strong regularity assumption}:
for some $C, \eta>0$, one has
{\renewcommand{\theequation}{\mbox{\texttt{SR}}}
\begin{equation}
\label{allmarginals}
\P \bigl\{ \|PX_i
\|_2^2 > t \bigr\} \le C t^{-1-\eta} \qquad\mbox{for } t >
C \rank(P)
\end{equation}}

\noindent for every orthogonal projection $P$ in $\R^n$.
Then, for $\varepsilon\in(0,1)$ and for
\[
N \ge C_\main \varepsilon^{-2-2/\eta}\cdot n,
\]
one has
%
%e1.4 #&#
\setcounter{equation}{3}
\begin{equation}
\label{maineqn} \E \Biggl\| \frac{1}{N} \sum_{i=1}^N
X_i X_i^T - I \Biggr\| \le\varepsilon.
\end{equation}
Here $C_\main= 512 (48C)^{2+2/\eta} (6+6/\eta)^{1+4/\eta}$, and as before
$\|\cdot\|$ denotes the spectral (operator) matrix norm, and $\|\cdot
\|_2$ denotes the Euclidean norm in~$\R^n$.
\end{theorem}

\begin{remark*}
Since the distribution of $PX_i$ is isotropic in the range of $P$, we have
$\E\|PX_i\|_2^2 = \rank(P)$.
This explains why (\ref{allmarginals}) concerns only
the tail values of $t$ which are above $\rank(P)$.
\end{remark*}

%s1.3 #&#
\subsection{Covariance estimation}

Returning to the covariance estimation problem, we deduce the following.

%co1.2 #&#
\begin{corollary}[(Covariance estimation)] \label{covarianceestimation}
Consider a random vector $X$ valued in $\R^n$ with covariance matrix
$\Sigma$.
Assume that for some $C,\eta>0$,
the isotropic random vector $Z = \Sigma^{-1/2} X$ satisfies
{\renewcommand{\theequation}{\mbox{\texttt{SR}}}
\begin{equation}
\label{allmarginalsnon-isotropic}
\P \bigl\{ \|PZ\|_2^2
> t \bigr\} \le C t^{-1-\eta} \qquad\mbox{for } t > C \rank(P)
\end{equation}
for every orthogonal projection $P$ in $\R^n$.
Then, for every $\varepsilon\in(0,1)$ and}
\[
N \ge C_\main \varepsilon^{-2-2/\eta}\cdot n,
\]
the sample covariance matrix $\Sigma_N$ obtained from $N$ independent
copies of $X$ satisfies
\[
\E\|\Sigma_N - \Sigma\| \le\varepsilon\|\Sigma\|.
\]
\end{corollary}

This result follows by applying Theorem~\ref{main} for the independent
copies of
the random vectors $Z_i= \Sigma^{-1/2} X_i$ instead of $X_i$, and by
multiplying the matrix
$\frac{1}{N} \sum_{i=1}^N X_i X_i^T - I$ in (\ref{maineqn}) by
$\Sigma^{1/2}$ on the left and
on the right.
Thus, for distributions satisfying (\ref{allmarginalsnon-isotropic})
we conclude that \textit{the minimal sample size for the covariance
estimation is
$N = O(n)$}.

Let us illustrate these results with two important examples.

%s1.4 #&#
\subsection{Sampling from log-concave distributions and convex sets}

A notable class of examples where Corollary~\ref{covarianceestimation} applies is
formed by the log-concave distributions, which includes the uniform
distributions
on convex bodies.
Consider a random vector $X$ with a log-concave distribution in $\R^n$, that is, whose density
has the form $e^{-V(x)}$ where $\log V(x)$ is a convex function on $\R^n$.
Paouris's concentration inequality~\cite{P} implies that
regularity assumption (\ref{allmarginalsnon-isotropic}) holds for $X$.
Indeed, consider an orthogonal projection $P$ in $\R^n$, and let $k =
\rank(P)$.
The distribution of the isotropic random vector $Z = \Sigma^{-1/2} X$
is log-concave in $\R^n$,
and so is the distribution of $PZ$ in the $k$-dimensional space $\range(P)$.
Paouris's theorem then states that
\[
\P \bigl\{ \|PZ\|_2^2 > t \bigr\} \le\exp(-ct)
\qquad\mbox{for
} t > Ck,
\]
where $C, c>0$ are absolute constants.
This is obviously stronger than assumption~(\ref{allmarginals}),
so Corollary~\ref{covarianceestimation} applies.

We conclude that \textit{the minimal sample size for estimating the
covariance matrix of a log-concave distribution is
$N = O(n)$}. This matches the bound obtained by Adamczak et al.~\cite{ALPT},
though it should be noted that the guarantee of~\cite{ALPT} holds with
probability that converges to $1$ exponentially
fast as $n \to\infty$, whereas ours holds only in expectation. We
have not tried to obtain probability bounds of this type;
note, however, that under our general assumption (\ref{allmarginals}),
the probability cannot converge to $1$ faster than at a polynomial
rate in $n$.

%s1.5 #&#
\subsection{Sampling from product distributions}

A distribution does not have to be log-concave in order to satisfy the
regularity assumptions
in Theorem~\ref{main} and Corollary~\ref{covarianceestimation}. For example,
all product distributions with finite $4+\varepsilon$ moments have
the required regularity property.
We can deduce this from the following thin shell estimate:

%pr1.3 #&#
\begin{proposition}[(Thin shell probability for product distributions)]
\label{productdistributions}
Let \mbox{$p \ge2$}, and consider a random vector $X = (\xi_1,\ldots,\xi_n)$,
where $\xi_i$ are independent random variables
with zero means, unit variances and with uniformly bounded $(2p)$th moments.
Then for every $1 \le k \le n$ and for every orthogonal projection $P$
in $\R^n$ with $\rank P = k$,
one has
%
%e1.5 #&#
\setcounter{equation}{4}
\begin{equation}
\label{productdistributionseqn} \E\bigl|\|PX\|_2^2 -
k\bigr|^p \lesssim k^{p/2}.
\end{equation}
The factor implicit in (\ref{productdistributionseqn}) depends only
on $p$ and on the bound on the $(2p)$th moments.
\end{proposition}

The proof of Proposition~\ref{productdistributions} is given in the
\hyperref[app]{Appendix}.

Applying Chebyshev's inequality together with (\ref{productdistributionseqn}),
we obtain for $t \ge k$ that
\[
\P \bigl\{ \|PX\|_2^2 > k + t \bigr\} \le
t^{-p} \cdot\E\bigl|\|PX\|_2^2 - k\bigr|^p
\lesssim t^{-p} k^{p/2} \le t^{-p/2}.
\]
Thus for $p>2$ we get a sub-linear tail, as required in the regularity
assumption (\ref{allmarginals}).

This shows that \textit{Theorem~\ref{main} applies for product
distributions in $\R^n$ with uniformly bounded $4+\varepsilon$ moments,
and it gives $N = O(n)$ for their covariance estimation}. Note that
this moment assumption is almost tight---according to~\cite{BSY}, if
the components $\xi_i$ are i.i.d. and have infinite fourth moment, then
$\limsup\|\Sigma_N\| \to\infty$ as $n \to\infty$ and $n/N \to y > 0$.
(This is because in this situation at least one of the $Nn$ i.i.d.
coordinates of $X_1, \ldots,X_N$ will likely to be large.)

%s1.6 #&#
\subsection{Extreme eigenvalues} \label{sextremeeigenvalues}

Theorem~\ref{main} states that, for sufficiently large~$N$, all
eigenvalues of the sample covariance matrix
$\Sigma_N = \frac{1}{N} \sum_{i=1}^N X_i X_i^T$ are concentrated
near $1$.
It is easy to extend this to a result that holds for all~$N$, as follows.

%co1.4 #&#
\begin{corollary} \label{fixedN}
Let $n, N$ be arbitrary positive integers,
suppose $X_i$ are independent isotropic random vectors in $\R^n$
satisfying (\ref{allmarginals}),
and let $y=n/N$. Then the sample covariance matrix
$\Sigma_N = \frac{1}{N}\sum_{i=1}^N X_iX_i^T$ satisfies
%
%e1.6 #&#
\begin{equation}
\label{newBai-Yin} 1 - C_1 y^c \le\E\l_{\min} (
\Sigma_N) \le\E\l_{\max} (\Sigma_N) \le1 +
C_1 \bigl(y+y^c\bigr).
\end{equation}
Here $c=\frac{\eta}{2\eta+2}$, $C_1=512 (16C)^{1+2/\eta} (6+6/\eta
)^{1+4/\eta}$
and $\l_{\min} (\Sigma_N)$, $\l_{\max} (\Sigma_N)$ denote the
smallest and
the largest eigenvalues of $\Sigma_N$, respectively.
\end{corollary}

We deduce this result in Section~\ref{supperedge}.
One can view (\ref{newBai-Yin}) as a nonasymptotic form of the
\textit{Bai--Yin law} for
the extreme eigenvalues of sample covariance matrices~\cite{BY}.
This law, associated with the works of
Geman, Bai, Yin, Krishnaiah and Silverstein
applies for product distributions, specifically for random vectors
$X = (\xi_1,\ldots, \xi_n)$ with i.i.d. components $\xi_i$ with
zero mean, unit variance
and finite fourth moment. For such distributions one has asymptotically
almost surely that
%
%e1.7 #&#
\begin{equation}
\label{Bai-Yin}
\bigl(1-\sqrt{y}\bigr)^2 - o(1) \le\l_{\min} (
\Sigma_N) \le\l_{\max} (\Sigma_N) \le\bigl(1+
\sqrt{y}\bigr)^2 + o(1)
\end{equation}
as $n \to\infty$ and $n/N \to y \in[0,1)$; see the rigorous
statement in~\cite{BY}.
This limit law is sharp. On the other hand, inequalities (\ref{newBai-Yin}) hold
in any fixed dimensions $N,n$ and for general distributions (as in
Theorem~\ref{main}),
\textit{without any independence requirements} for the coordinates.

\begin{remark*}
Comparing (\ref{newBai-Yin}) with (\ref{Bai-Yin}) one can ask about
the optimal
value of the exponent $c$, in particular whether $c=1/2$.
In a recent paper~\cite{ALPTsharp}, Adamczak et al.
obtained the optimal exponent $c=1/2$ for log-concave distributions,
and more generally for sub-exponential distributions in the sense of
(\ref{sub-exponential}). As (\ref{sub-exponential}) implies (\ref
{allmarginals})
with $\eta=(p-1)/2$ and $C \le(O(p))^p$, Theorem~\ref{main} recovers
a bound of
$c=1/2-1/(p+1)=1/2-o(1)$ as $p\to\infty$.
\end{remark*}

\begin{remark*}[(Random matrices with independent rows)]
Corollary~\ref{fixedN} can be interpreted as a result about the
spectrum of
random matrices with independent rows.
Indeed, if $A$ is the matrix with rows $X_i$, then
$\Sigma_N = \frac{1}{N} \sum_{i=1}^N X_i X_i^T = \frac{1}{N} A^T A$.
So the singular values of the matrix $\frac{1}{\sqrt{N}} A$
are the same as the eigenvalues of the matrix $\Sigma_N$,
and they are controlled as in (\ref{newBai-Yin}). In particular,
under the regularity assumption (\ref{allmarginals}) on $X_i$ we
obtain that
\[
\bigl(\E\|A\|^2\bigr)^{1/2} \le C_2 \bigl(\sqrt{N} +
\sqrt{n}\bigr),
\]
where $C_2 = \sqrt{2C_1}$, and $C_1$ is as in Corollary~\ref{fixedN}.
%Let us illustrate this point with an optimal estimate on the spectral
%norm of $A$.
% Let $n, N$ be arbitrary positive integers,
% and suppose $X_i$ are independent isotropic random vectors in $\R^n$
%satisfying \eqref{allmarginals}.
% Consider the $N \times n$ random matrix $A$ with rows $X_1, \ldots,
%X_N$.
% Then
% $$
% (\E\|A\|^2)^{1/2} \le C_2 (\sqrt{N} + \sqrt{n})
% $$
% where $C_2 = \sqrt{2C_1}$ and $C_1$ is as in Corollary~\ref{fixedN}.
%The result follows from the identity
%$\l_{\max}(\Sigma_N) = \|\Sigma_N\| = \frac{1}{N} \|A\|^2$
%and the upper bound in \eqref{newBai-Yin}.

Notice that while the rows of matrix $A$ are independent, \textit{the
columns of $A$ may be dependent}. The simpler case where all
\textit{entries} of $A$ are independent is well understood by now. In the
latter case, if the entries have zero mean and uniformly bounded fourth
moments, the bound $\E\|A\| \lesssim\sqrt{N} + \sqrt{n}$ follows, for
example, from Latala's general inequality~\cite{L}.
%Furthermore, Corollary~\ref{norms} was obtained in~\cite{ALPT} in the
%special case
%where the rows $X_i$ have log-concave distributions and, more
%generally, if $X_i$ satisfy
%the sub-exponential moment assumptions \eqref{sub-exponential}.
\end{remark*}

%s1.7 #&#
\subsection{Smallest eigenvalue}

Our proof of Theorem~\ref{main} consists of two separate arguments for upper
and lower bounds for the spectrum of the sample covariance matrix.
It turns out that the full power of the strong regularity assumption
(\ref{allmarginals}) is not needed for the lower bound.
It suffices to assume \textit{$2+\eta$ moments
for one-dimensional marginals} rather than for marginals in all dimensions.
This is only slightly stronger than the isotropy assumption, which
fixes the
\textit{second moments} of one-dimensional marginals, and it broadens the class
of distributions for which the result applies.
We state this as a separate theorem.

%th1.5 #&#
\begin{theorem}[(Smallest eigenvalue)] \label{mainlower}
Consider independent isotropic random vectors $X_i$ valued in $\R^n$.
Assume that $X_i$ satisfy the following \textup{weak regularity
assumption}: for some
$C,\eta>0$,
{\renewcommand{\theequation}{\mbox{\texttt{WR}}}
\begin{equation}
\label{onedimmarginals}
\sup_{\|x\|_2 \le1} \E\bigl|\langle
X_i, x\rangle\bigr|^{2+\eta} \le C.
\end{equation}
Then, for $\varepsilon> 0$ and for}
%
%e1.8 #&#
\setcounter{equation}{7}
\begin{equation}
\label{Cprime2} N \ge C_\low \varepsilon^{-2-2/\eta} \cdot n,
\end{equation}
the minimum eigenvalue of the sample covariance matrix $\Sigma_N =
\frac{1}{N}\sum_{i=1}^N X_iX_i^T$
satisfies
\[
\E\l_{\min} (\Sigma_N) \ge1 - \varepsilon.
\]
Here $C_\low= 40(10C)^{2/\eta}$.
\end{theorem}

%equivalent to the following statement:
% for all $N \ge1$, we have
% $$
% \E\l_{\min} (\Sigma_N) \ge1 - \Big(C_\low\frac{n}{N}\Big)^\frac{
% $$

% Note that the isotropy condition $\E X_i X_i^T = I$ alone implies
%that \eqref{onedimmarginals}
% holds for $\eta= 0$:
%$$
%= \sup_{\|x\|_2 \le1} \|x\|^2 = 1.
%$$
%So regularity assumption \eqref{onedimmarginals} is a bit stronger
%than the
%trivial variance identity.

\begin{remark*}[(Moments vs. tails)]
We have chosen to write (\ref{onedimmarginals}) in terms of
moments rather than in terms of tail bounds as in (\ref{allmarginals}).
By integration of the tails one can check that, for any given $\eta>0$,
(\ref{allmarginals}) with parameter $C$ implies (\ref{onedimmarginals})
with parameter $C'=C(2+2/\eta)$.

In the remainder of the paper we will use (\ref{onedimmarginals})
for theorems
regarding only the smallest eigenvalue and (\ref{allmarginals}) for theorems
which involve the largest one.
\end{remark*}

\begin{remark*}[(Product distributions with $2+\eta$ moments)]
Many distributions of interest satisfy (\ref{onedimmarginals}).
For example, let $X = (\xi_1,\ldots, \xi_n)$ have i.i.d. components
$\xi_i$
with zero mean, unit variance and finite $(2+\eta)$ moment.
Then a standard application of symmetrization and Khintchine's inequality
(or a direct application of Rosenthal's inequality~\cite{Ros}, see
\cite{FHJSZ})
shows that one-dimensional
marginals of $X$ also have bounded $(2+\eta)$ moments; that is,
(\ref{onedimmarginals}) holds.

In the context of the Bai--Yin law discussed in Section~\ref{sextremeeigenvalues},
this indicates that the \textit{smallest} eigenvalue of a random matrix
can be approximately controlled [as in (\ref{newBai-Yin})]
even if the \textit{fourth moment is infinite}. However, as we already
recalled, four moments are necessary
to control the \textit{largest} eigenvalue in the classical Bai--Yin law
\cite{BSY}.
\end{remark*}

\begin{remark*}[(Covariance estimation)]
Theorem~\ref{mainlower} can be used to obtain a \textit{lower} estimate
for the covariance matrix
under the weak regularity assumption (\ref{onedimmarginals}).

%This the following is a one-sided
%version of Corollary~\ref{covarianceestimation}.
%Consider a random vector $X$ valued in $\R^n$ with covariance matrix $
%Assume that for some $C,\eta>0$, the isotropic random vector $Z =
%$$
%$$
%Then, for $\e>0$, the sample covariance matrix $\Sigma_N$ obtained
%from $N\ge C_\low n/\e^{2+2/\eta}$
%independent copies of $X$ satisfies
%$$
%$$ where $\succeq$ denotes the canonical order ($A \succeq B$ if $A-B$
%is positive semidefinite).
\end{remark*}

%s1.8 #&#
\subsection{Optimality of the regularity assumptions} \label{soptimality}

Let us briefly mention two simple and known examples that illustrate
the role of regularity assumptions
(\ref{allmarginals}) and (\ref{onedimmarginals}) in the control of
the largest and
smallest eigenvalues, respectively.

For the largest eigenvalue as in Theorem~\ref{main}, it is not
sufficient to put a regularity
assumption of the type (\ref{allmarginals}) only on
\textit{one-dimensional} marginals,
as it is done in Theorem~\ref{mainlower} for the smallest eigenvalue.
Even the following very strong (exponential) moment assumption is insufficient:
%
%e1.9 #&#
\begin{equation}
\label{exponential} \sup_{\|x\|_2 \le1} \P \bigl\{ \bigl|\langle X, x\rangle\bigr| > t
\bigr\} \le C \exp(-ct) \qquad\mbox{for } t > 0.
\end{equation}
Indeed, consider a random vector $X = \xi Z$ where
$Z$ is a random vector uniformly distributed in the Euclidean sphere in
$\R^n$ centered at the origin
and with radius~$\sqrt{n}$, and where $\xi$ is a standard normal
random variable.
Then $X$ is isotropic, and all one-dimensional marginals of $X$ have exponential
tail decay (\ref{exponential}). However, the multiplier $\xi$
produces a dimension-free
tail decay of the norm of $Z$, namely
$\P\{ \|X\|_2 > t \sqrt{n} \} = \P\{ \xi> t \}
\gtrsim\exp(-C't^2)$ for $t > 0$.
It follows that a sample of $N$ independent copies $X_1, \ldots, X_N$
of $X$
safisfies $\E\max_{i \le N} \|X_i\|_2^2 \gtrsim N \log N$, so the
matrix $\Sigma_N = \frac{1}{N} \sum_{i=1}^N X_i X_i^T$
satisfies
\[
\E\|\Sigma_N - I\| \ge N^{-1} \E\max_{i \le N}
\|X_i\|_2^2 - 1 \gtrsim\log N,
\]
which contradicts the conclusion of Theorem~\ref{main}.
This example is essentially due to Aubrun; see~\cite{ALPT},
Remark 4.9.

\begin{remark*}
It is not clear whether Theorem~\ref{main} would hold if, in addition
to $(2+\eta)$ moments on one-dimensional marginals, one puts a total
boundedness assumption
\[
\|X\| = O\bigl(\sqrt{n}\bigr) \qquad\mbox{almost surely.}
\]
A conjecture of this type is discussed in~\cite{Vcovariance} where a
version of
the theorem is proved under this assumption, with $\eta=2$ but with an
additional
$(\log\log n)^{O(1)}$ oversampling factor.
\end{remark*}

Furthermore, we note that for the smallest eigenvalue as in
Theorem~\ref{mainlower}, one cannot drop
the regularity assumption (\ref{onedimmarginals}); that is, the
assumption with $\eta= 0$ is not sufficient.
This is seen for $X_i$ uniformly distributed in the set of $2n$ points
$(\pm e_k)$ where $(e_k)_{k=1}^n$ is
an orthonormal basis in $\R^n$. Indeed, in order that the smallest
eigenvalue of the matrix
$\Sigma_N = \frac{1}{N} \sum_{i=1}^N X_i X_i^T$ be different from
zero, one needs $\Sigma_N$ to have
full rank, for which all $n$ basis vectors $e_k$ need be present in the
sample $X_1, \ldots, X_N$.
By the coupon collector's problem, for this to happen with constant
probability one needs
a sample of size $N \gtrsim n \log n$.
For $N = o(n \log n)$, the smallest eigenvalue is zero with high
probability, so the conclusion
of Theorem~\ref{mainlower} fails.

%s1.9 #&#
\subsection{The argument: Randomizing the spectral sparsifier}
\label{secargsketch}

Our proof of Theorem~\ref{main} consists of randomizing the
\textit{spectral sparsifier}
invented by Batson, Spielman and Srivastava~\cite{BSS}; see
\cite{SPhD}.
The randomization makes the spectral sparsifier appear naturally
in the context of random matrix theory. The method is based on
evaluating the Stieltjes transform of $\Sigma_N$ while making rank one updates.
However, in contrast to typical methods of random matrix theory (and to
the spectral sparsifier itself), we shall
evaluate the Stieltjes transform \textit{at random real points}.

Let us illustrate the method by working out a crude upper bound $O(1)$
for the largest eigenvalue of $\Sigma_N$. Equivalently, we want to
show that a
general Wishart matrix $A_N:= N \Sigma_N = \sum_{i=1}^N X_i X_i^T$
has all eigenvalues bounded by $O(N)$.
We evaluate the Stieltjes transform
%
%e1.10 #&#
\begin{equation}
\label{Stieltjes} m_{A_N}(u) = \tr(uI-A_N)^{-1} =
\sum_{i=1}^n \bigl( u -
\l_i(A_N) \bigr)^{-1},\qquad u \in\R,
\end{equation}
where $\l_i(A_N)$ denote the eigenvalues of $A_N$. This function has
singularities at the points $\l_i(A_N)$, and it
vanishes at infinity. So the largest eigenvalue of $A_N$ is the largest
$u$ where $m_{A_N}(u) = \infty$.
However, such $u$ is difficult to compute. So we soften this quantity
by considering the
largest number $u_N$ that satisfies
%
%e1.11 #&#
\begin{equation}
\label{eqspectraledge} m_{A_N}(u_N) = \pn,
\end{equation}
where $\pn$ is a fixed sensitivity parameter, for example, $\pn=1$.
%At this point a physical interpretation can be helpful. Let us put
%unit electric charges at the points $\l_i(A)$.
%Then the Stieltjes transform $m_{A_N}(u)$ is the electric {\em
%potential} evaluated at the point $u$.
%The point $u_N$ is the rightmost point with a given reading $\pn$. The
%larger the sensitivity value $\pn$,
%the closer the point $u_N$ is to the largest eigenvalue $\lmax(A_N)$
%that we seek to control; the smaller $\pn$, the softer and
%easier to control $u_N$ becomes.

The \textit{soft spectral edge} $u_N$ provides an upper bound for the
actual spectral edge, $\lmax(A_N) < u_N$.
So our goal is to show that
\[
\E u_N = O(N).
\]
This is the same problem as in~\cite{BSS}, except the eigenvalues and hence
the soft spectral edge $u_N$ are now \textit{random} points. The
randomized problem is more difficult as we note below.

As opposed to the largest eigenvalue of $A$, the soft spectral edge
$u_N$ can be computed inductively using rank-one updates to the matrix;
$u_N$ will move to the right by a random amount at each step
as we replace $A_{k-1}$ by $A_k = A_{k-1} + X_k X_k^T$.
Initially, $A_0$ = 0 so $u_0=n$.
It suffices to prove that the $u_k$ moves by $O(1)$ on average at each
step:
%
%e1.12 #&#
\begin{equation}
\label{barriershiftsmall} \E(u_k - u_{k-1}) = O(1).
\end{equation}
Indeed, by summing up we would obtain the desired estimate
$\E u_N = n + O(1)N = O(N)$.

The soft edge $u_k$ can be recomputed at each step because it is
determined by the Stieltjes transform $m_{A_k}(u)$,
which in turn can be recomputed using Sherman--Morrison
formula,
as is done in~\cite{BSS}, which gives for every $u \in\R$ that
%
%e1.13 #&#
\begin{equation}
\label{eqSherman-Morrison} m_{A_k}(u) = m_{A_{k-1}}(u) +
\frac{X_k^T (u I - A)^{-2} X_k}{1 -
X_k^T (u I - A)^{-1} X_k}.
\end{equation}
This reduces proving (\ref{barriershiftsmall}) to a probabilistic
problem, which is essentially governed by the distribution
of the random vector $X_k$.

The difficulty is that we are facing a nonlinear inverse problem.
Indeed, for a fixed $u$ it is not difficult to compute
the expectation of $m_{A_k} (u)$ from (\ref{eqSherman-Morrison}), and
in particular to bound the expectation by $\phi$;
this is done in~\cite{BSS}.
However, we require the identity $m_{A_k}(u) = \phi$ to hold
\textit{deterministically},
because the largest $u$ that satisfies it defines the soft spectral
edge of $A_k$
as in (\ref{eqspectraledge}).
The task of computing the expectation of a random number $u$ for which
$m_{A_k}(u) = \phi$
is a highly nonlinear inverse problem~\cite{BN}, Section 4.1.
This is where some regularity of $X_k$ with respect to the
eigenstructure of $A_{k-1}$ becomes essential.
A technical part of our argument developed in most of the remaining sections
is to realize and prove that a small amount or regularity encoded by
(\ref{allmarginals}) or (\ref{onedimmarginals}) is already
sufficient to control the solution to the inverse problem,
and ultimately to control the spectral edges of $A$.

%s1.10 #&#
\subsection{Organization of the paper}

The rest of the paper is organized as follows. We start with the
somewhat simpler Theorem~\ref{mainlower}
for the smallest eigenvalue in Section~\ref{sloweredge}. A
corresponding result for the
largest eigenvalue, Theorem~\ref{mainupper}, is proved in
Section~\ref{supperedge}.
Corollary~\ref{fixedN} is also deduced in Section~\ref{supperedge}.
Combining Theorems~\ref{mainlower} and~\ref{mainupper} in
Section~\ref{sspectralnorm},
we obtain the main Theorem~\ref{main} on the spectral norm. In the
\hyperref[app]{Appendix}, we prove Proposition~\ref{productdistributions} on the
regularity of product distributions.

%s2 #&#
\section{The lower edge} \label{sloweredge}

We begin by proving Theorem~\ref{mainlower} about the the lower edge of
the spectrum, which is slightly simpler and requires fewer assumptions
than the
upper edge. As in~\cite{BSS}, the tool that we use to do this is the
\textit{lower Stieltjes transform}
\[
\mn_A(\ell) = \tr(A-\ell I)^{-1} = \sum
_{i=1}^n \bigl( \lambda_i(A)-\ell
\bigr)^{-1},\qquad \ell\in\R.
\]
Note that $\mn_A(\ell) = -m_{-A}(-\ell)$ where $m_A$ is
the usual Stieltjes transform in~(\ref{Stieltjes}).

For a sensitivity value $\pn> 0$, we define the
\textit{lower soft spectral edge}
$\ell_\pn(A)$ to be the smallest $\ell$ for which
\[
\mn_A(\ell) = \pn.
\]
Since $\mn_A(\ell)$ increases from $0$ to $\infty$ as $\ell$
increases from $-\infty$ to the lower
spectral edge $\lmin(A)$, the value $\ell_\pn(A)$ is defined
uniquely, and we always have the bound
\[
\ell_\pn(A) < \lmin(A).
\]
For $\pn\to\infty$ we have $\ell_\pn(A) \to\lmin(A)$. However,
we will work with
small sensitivity $\pn\in(0,1)$, which will make the soft spectral
edge $\ell_\pn(A)$ softer and easier to control.

The crucial property of $\ell_\pn(A)$ is that it grows steadily under
rank-one updates.
Consider what happens when we add a random rank-one matrix $XX^T$ to
$A\succ\ell I$, where $X$ is chosen from an isotropic distribution on
$\R^n$.
As $\E\tr(A+XX^T)=\tr(A)+\tr\E XX^T = \tr(A)+n$,
we expect the eigenvalues of
$A+XX^T$ to have increased by $1$ on average.
It turns out that $\ell_\pn(A)$ behaves almost as nicely as this
if the distribution of $X$ is sufficiently regular and the
sensitivity $\pn$ is sufficiently small. This is established in the
following theorem.

%th2.1 #&#
\begin{theorem}[(Random lower shift)] \label{thmlowerrankone}
Suppose
$X$ is an isotropic random vector in $\R^n$ satisfying the weak
regularity assumption:
for some $C,\eta>0$,
{\renewcommand{\theequation}{\mbox{\texttt{WR}}}
\begin{equation}
\label{eonedimmarginals}
\sup_{\|x\|\le1} \E\bigl|\langle X,x
\rangle\bigr|^{2+\eta}\le C.
\end{equation}
Let $\varepsilon> 0$ and}
\[
\pn\le c_{\rref{thmlowerrankone}} \varepsilon^{1+2/\eta},
\]
where $c_{\rref{thmlowerrankone}}^{-1} = 10(5C)^{2/\eta}$.
Then for every symmetric $n \times n$ matrix $A$, one has
\[
\E\ell_\pn\bigl(A+XX^T\bigr) \ge\ell_\pn(A)
+ 1-\varepsilon.
\]
\end{theorem}

Iterating Theorem~\ref{thmlowerrankone} easily yields a proof of
Theorem~\ref{mainlower} as
follows.
\begin{pf*}{Proof of Theorem~\ref{mainlower}}
Let $A_0=0$ and $A_k=A_{k-1}+X_kX_k^T$ for \mbox{$k\le N$}.
Setting $\pn= c_{\rref{thmlowerrankone}}   \varepsilon^{1+2/\eta
}$, we find that
\[
\ell_\pn(A_0) = \frac{-n}{\pn}. %= -n\frac{10(5C)^{2/\eta}}{\e^{1+2/
\]
Applying Theorem~\ref{thmlowerrankone} inductively to $A_0,
A_1,\ldots,A_N$, we find that
\[
\E \bigl[ \ell_\pn(A_k)-\ell_\pn(A_{k-1})
| A_{k-1} \bigr] \ge 1-\varepsilon \qquad\mbox{ for all $k\le N$},
\]
where we take the conditional expectation with respect to the random
vector~$X_k$, given the random vectors
$X_1,\ldots,X_{k-1}$, that is, given $A_{k-1}$.
Summing up these bounds yields
%
%e2.1 #&#
\setcounter{equation}{0}
\begin{equation}
\label{lminAN} \E\ell_\pn(A_N) \ge\ell_\pn(A_0)
+N(1-\varepsilon).
\end{equation}
Recalling that $\lmin(A_N) > \ell_\pn(A_N)$ and dividing both sides
of (\ref{lminAN}) by $N$, we conclude that
\[
\E\lmin \Biggl( \frac{1}{N}\sum_{i=1}^N
X_iX_i^T \Biggr) > \frac{\ell_\pn(A_0)}{N} + 1-
\varepsilon= 1 - \varepsilon- \frac{n}{\pn N}.
\]
For $N \ge n/\varepsilon\pn$, the bound becomes $1-2\varepsilon$.
Substituting the value of $\pn$ and replacing $\varepsilon$ by
$\varepsilon/2$ gives the promised result.
\end{pf*}

%independence of the random vectors $X_i$.
% To invoke Theorem~\ref{thmlowerrankone}, we only need that the
%distribution
% of each $X_i$ conditioned on $X_1,\ldots, X_{i-1}$ satisfies (

The rest of this section is devoted to proving Theorem~\ref{thmlowerrankone}.
Given a matrix~$A$, a~real number $\ell< \lmin(A)$ and a vector $x
\in\R^n$,
we say that $\dn\ge0$ is a \textit{feasible lower shift} if
\[
A\succ(\ell+\dn)I \quad\mbox{and}\quad \mn_{A+xx^T}(\ell+\dn) \le
\mn_A(\ell).
\]
The definition of the soft spectral edge $\ell= \ell_\pn(A)$ along
with monotonicity of the Stieltjes transform implies
that
\[
\ell_\pn\bigl(A+xx^T\bigr)\ge\ell_\pn(A)+\dn
\]
for every feasible lower shift $\dn$. So we will be done if we can
produce a feasible shift $\dn$
such that $\E\dn\ge1-\varepsilon$ where the expectation is over
random $X$.

We begin by reducing the feasibility for a shift $\dn$ to an
inequality involving two quadratic forms. The following lemma
appeared in~\cite{BSS}, and we include it with a proof for completeness.

%le2.2 #&#
\begin{lemma}[(Feasible lower shift)] \label{lowerfeasible}
Consider the numbers $\ell\in\R$, $\delta> 0$, a~matrix $A\succ
(\ell+\dn) I$ and a vector $x$.
Then a sufficient condition for
%
%e2.2 #&#
\begin{equation}
\label{mlem} \mn_{A+xx^T}(\ell+\dn) \le\mn_A(\ell)
\end{equation}
is\setcounter{footnote}{2}\footnote{To ease the notation, we sometimes write $A-u$ instead of $A-uI$.}
%
%e2.3 #&#
\begin{equation}
\label{lowerbarriercondition}
\qquad\frac{1}{\dn}\frac{x^T(A-\ell-\dn)^{-2}x}{\tr(A-\ell-\dn)^{-2}} -
x^T(A-\ell-\dn)^{-1}x =: \frac{1}{\dn}
\Qn_2(\dn,x)- \Qn_1(\dn,x) \ge1.
\end{equation}
\end{lemma}

\begin{pf}
We begin by expanding $\mn_{A+xx^T}(\ell+\dn)$ using the
Sherman--Morisson formula,
\begin{eqnarray*}
\mn_{A+xx^T}(\ell+\dn) &=& \tr\bigl(A+xx^T-\ell-\dn
\bigr)^{-1} \\
&=&\tr(A-\ell-\dn)^{-1} - \frac{x^T(A-\ell-\dn)^{-2}x}{1+x^T(A-\ell
-\dn)^{-1}x}.
\end{eqnarray*}
Furthermore,
\[
\tr(A-\ell-\dn)^{-1} = \mn_A(\ell) + \tr\bigl[(A-\ell-
\dn)^{-1}-(A-\ell)^{-1}\bigr].
\]
The assumption $A\succ(\ell+\dn)I$ implies that
\[
(A-\ell-\dn)^{-1}-(A-\ell)^{-1}\preceq\dn(A-\ell-
\dn)^{-2}.
\]
Combining these estimates, we see that (\ref{mlem}) holds as long as
\[
\delta\cdot\tr(A-\ell-\dn)^{-2} - \frac{x^T(A-\ell-\dn
)^{-2}x}{1+x^T(A-\ell-\dn)^{-1}x} \le0,
\]
which we can rearrange into (\ref{lowerbarriercondition}) observing
that all quadratic
forms involved are positive.
\end{pf}

Inequality (\ref{lowerbarriercondition}) is quite nontrivial in the
sense that $\dn$ appears in many places, and it is not immediately
clear from looking at it what
the largest feasible $\dn$ is given $A,x$ and $\ell$. In the
following lemma,
we present a tractable and explicit quantity defined solely in terms of
$\Qn_1(0,x)$ and $\Qn_2(0,x)$ which
always satisfies (\ref{lowerbarriercondition}) and thus provides a
lower bound on the best possible~$\dn$.

%le2.3 #&#
\begin{lemma}[(Explicit feasible shift)]\label{dnexplicit}
Consider numbers $\ell\in\R$, $\pn> 0$, a~matrix $A\succ\ell I$
satisfying $\mn_A(\ell)\le\pn$,
and a vector $x$.
Then for every $t \in(0,1)$, the shift
\[
\dn:= (1-t)^3\Qn_2(0,x) \one_{\{ \Qn_1(0,x)\le t\}}
\one_{\{
\Qn_2(0,x)\le t/\pn\}}
\]
satisfies $A \succ(\ell+ \dn) I$ and condition (\ref
{lowerbarriercondition}).
Therefore $\dn$ is a feasible lower shift, that is, $\mn_{A+xx^T}(\ell+\dn)\le\mn_A(\ell)$.
\end{lemma}

The proof is based on regularity properties of the quadratic forms
$q_1$ and~$q_2$, which we state in the following
two lemmas.

%le2.4 #&#
\begin{lemma}[(Regularity of quadratic forms)]
\label{lemlowerquadraticregularity}
Consider the numbers \mbox{$\ell\in\R$}, $\pn> 0$, a matrix $A\succ\ell
I$ satisfying $\mn_A(\ell)\le\pn$,
and a vector $x$. Then for every positive number $\delta< 1/\pn$,
one has
$A \succ(\ell+\dn)I$, and moreover:
\begin{longlist}
\item$\Qn_1(0,x) \le\Qn_1(\dn,x) \le(1-\dn\pn)^{-1}\Qn_1(0,x)$;
\item$(1-\dn\pn)^2 \Qn_2(0,x) \le\Qn_2(\dn, x) \le(1-\dn\pn
)^{-2}\Qn(0,x)$.
\end{longlist}
\end{lemma}

\begin{pf}
The assumption $A \succ\ell I$ states that all eigenvalues $\l_i$ of $A$
satisfy $\l_i > \ell$.
Together with the assumption $\mn_A(\ell) =\sum_i (\lambda_i-\ell
)^{-1} \le\pn$ this implies
that $(\lambda_i-\ell)^{-1} \le\pn$ for all $i$, and hence $\l_i -
\ell\ge1/\pn>
\dn$ and $A \succ(\ell+\dn)I$ as claimed.

\mbox{}\hphantom{i}(i)
Let $(\psi_i)_{i \le n}$ denote the eigenvectors of $A$; then
%
%e2.4 #&#
\begin{equation}
\label{q1diag} \Qn_1(\delta,x) = \sum_{i=1}^n
\frac{\langle x,\psi_i\rangle^2}{\lambda_i-\ell-\dn}.
\end{equation}
Recalling that $\l_i - \ell\ge1/\pn$, we have the comparison inequalities
\[
(1-\dn\pn) (\lambda_i-\ell) =\lambda_i-\ell-\pn\dn(
\lambda_i-\ell) \le\lambda_i-\ell-\dn \le
\lambda_i-\ell.
\]
Using these for every term in (\ref{q1diag}), we complete the proof
of (i).\vadjust{\goodbreak}

(ii) Similar to (i), noting that the numerator and
denominator of $\Qn_2$ are increasing in $\delta$.
\end{pf}

%le2.5 #&#
\begin{lemma}[(Moments of quadratic forms)] \label{momentsquadratic}
Consider numbers $\ell\in\R$, $\pn> 0$ and a matrix $A\succ\ell I$
satisfying $\mn_A(\ell) \le\pn$.
If $X$ is an isotropic random vector satisfying (\ref
{eonedimmarginals}), then for $p=1+\eta/2$ the following
moment bounds hold:
\begin{longlist}
\item$\E\Qn_1(0,X) = \mn_A(\ell) \le\pn$ and $\E\Qn_1(0,X)^p
\le C\pn^p$;
\item$\E\Qn_2(0,X) = 1$ and $\E\Qn_2(0,X)^p \le C$.
\end{longlist}
\end{lemma}

\begin{pf}
(i) As in the proof of the previous lemma, let $(\psi_i)_{i \le
n}$ denote the eigenvectors of $A$.
By isotropy we have
\[
\E\Qn_1(0,X) = \sum_{i=1}^n
\frac{\E\langle
X,\psi_i\rangle^2}{\lambda_i-\ell} = \mn_A(\ell) \le\pn.
\]
For the moment bound we use Minkowski's inequality to obtain
\[
\bigl(\E\Qn_1(0,X)^p\bigr)^{1/p} \le\sum
_{i=1}^n \frac{(\E\langle X,\psi_i\rangle^{2p})^{1/p}}{\lambda_i-\ell} \le\sum
_{i=1}^n \frac{C^{1/p}}{\lambda_i-\ell} =
C^{1/p} \mn_A(\ell) \le C^{1/p} \pn.
\]

(ii) Analogous to (i).
\end{pf}

We can now finish the proof of Lemma~\ref{dnexplicit}.

\begin{pf*}{Proof of Lemma~\ref{dnexplicit}}
First observe that by construction,
%
%e2.5 #&#
\begin{equation}
\label{deltasmall} \dn\le\Qn_2(0,x) \one_{\{ \Qn_2(0,x)\le t/\pn\}} \le t/\pn<
1/\pn,
\end{equation}
so that we always have $A\succ(\ell+\dn)I$ by Lemma
\ref{lemlowerquadraticregularity}.

If either of the indicators in the definition of the shift $\delta$
is zero, then
$\delta=0$, which is trivially feasible, and we are done. So assume
both indicators are nonzero, that is,
$\Qn_1(0,x)\le t$ and $\Qn_2(0,x)\le t/\pn$. By Lemma
\ref{lowerfeasible}, it suffices to prove
inequality~(\ref{lowerbarriercondition}), which is equivalent to
\[
\frac{\Qn_2(\dn,x)}{1+\Qn_1(\dn,x)} \ge\delta.
\]
We can show this by replacing $\dn$ with zero using Lemma
\ref{lemlowerquadraticregularity}:
\begin{eqnarray*}
\frac{\Qn_2(\dn,x)}{1+\Qn_1(\dn,x)} &\ge&\frac{\Qn_2(0,x)(1-\dn\pn)^2}{1+\Qn_1(0,x)(1-\dn\pn)^{-1}}
\\
&\ge&\frac{\Qn_2(0,x)(1-t)^2}{1+t(1-t)^{-1}}\qquad
\mbox{[as $\dn\pn\le t$ by (\ref{deltasmall}) and $
\Qn_1(0,x)\le t$]}
\\
&=& \Qn_2(0,x) (1-t)^3 = \delta.
\end{eqnarray*}

The proof is complete.\vadjust{\goodbreak}
\end{pf*}

We now complete the proof of Theorem~\ref{thmlowerrankone}
by using the regularity properties of $X$ to show that the
expectation of $\dn$, as defined in Lemma~\ref{dnexplicit}, is large.
Roughly speaking, this happens because (1) $\dn$ is defined to be
slightly less than
$\Qn_2(0,X)$ whenever \textit{both} $\Qn_1(0,X)$ and $\Qn_2(0,X)$ are
not too large;
(2) that event occurs with very high probability when
$\pn$ is sufficiently small; (3) the expectation of $\Qn_2(0,X)$
equals $1$.

\begin{pf*}{Proof of Theorem~\ref{thmlowerrankone}}
Let $\ell= \ell_\pn(A)$; then $\mn_A(\ell) = \pn\le c_{\rref
{thmlowerrankone}} \varepsilon^{1+2/\eta}$ by assumption.
Define a feasible shift $\delta$ as in Lemma~\ref{dnexplicit} for $t
= \varepsilon/5$.
Recall that it suffices to prove that $\E\dn\ge1-\varepsilon$.

According to Lemma~\ref{dnexplicit},
\begin{eqnarray*}
\E\dn&=& {(1-t)^3} \bigl[ \E\Qn_2(0,X) - \E
\Qn_2(0,X) \one_{\{
\Qn_1(0,X)> t \lor\Qn_2(0,X) > t/\pn\}} \bigr]
\\
& \ge&(1-t)^3 \bigl[ 1 - \bigl(\E\Qn_2(0,X)^{p}
\bigr)^{1/p} \cdot\bigl(\P {\bigl\{\Qn_1(0,X)> t \lor
\Qn_2(0,X) > t/\pn\bigr\}}\bigr)^{1/q} \bigr],
\end{eqnarray*}
where we used H\"older's inequality with exponents $p=1+\eta/2$ and
$q=\frac{p}{p-1} =2/\eta+1$.
By Lemma~\ref{momentsquadratic}, we have $\E\Qn_2(0,X)^{p} \le C$.
Next, the probability can be estimated by a union bound, Markov's
inequality and the moment bounds
of Lemma~\ref{momentsquadratic}, which gives
\begin{eqnarray*}
&&\P\bigl\{\Qn_1(0,X)> t \lor \Qn_2(0,X) > t/\pn\bigr\}
\\
&&\qquad\le\P\bigl\{\Qn_1(0,X)^p > t^p\bigr\} + \P
\bigl\{\Qn_2(0,X)^p > (t/\pn)^p\bigr\}
\\
&&\qquad\le\frac{C\pn^p}{t^p} + \frac{C}{(t/\pn)^p} = 2C (\pn/t)^p.
\end{eqnarray*}
We conclude that
\begin{eqnarray*}
\E\delta &\ge& (1-t)^3 \bigl[ 1 - C^{1/p} \cdot\bigl(2C (
\pn/t)^p\bigr)^{1/q} \bigr]
\\
&\ge& (1-t)^3 \bigl[ 1 - 2C (\pn/t)^{\eta/2} \bigr]
\qquad\mbox{(as $1/p+1/q=1$ and $p/q=\eta/2$)}
\\
&=& (1-\varepsilon/5)^3 (1-\varepsilon/5) \qquad\mbox{(substituting $t$ and
the bound for $\pn$)}
\\
&\ge& 1-\varepsilon
\end{eqnarray*}
as promised.
\end{pf*}

%It is easy to check that the above proof of Theorem
%an {\em approximate} isotropy condition
%I,\end{equation}
%for some constants $c_1,c_2$ depending on $\e$.
%Thus if we desire a bound of $\lmin(\Sigma_N) >
%1-\e$ in Theorem~\ref{mainlower}, then $\E X_iX_i^T = I$ can be
%replaced by the
%weaker condition \eqref{weakisotropy}.

%s3 #&#
\section{The upper edge} \label{supperedge}

In this section we establish the following estimate for the expected
largest eigenvalue,
analogous to Theorem~\ref{mainlower} for the smallest one.

%th3.1 #&#
\begin{theorem}[(Largest eigenvalue)] \label{mainupper}
Consider independent isotropic random vectors $X_i$ valued in $\R^n$.
Assume that $X_i$ satisfy (\ref{allmarginals}) for some $C,\eta>0$.
Then, for $\varepsilon\in(0,1)$ and for
\[
N \ge C_\upper \varepsilon^{-2-2/\eta} \cdot n,
\]
the maximum eigenvalue of the sample covariance matrix $\Sigma_N =
\frac{1}{N} \sum_{i=1}^N X_i X_i^T$ satisfies
%
%e3.1 #&#
\begin{equation}
\label{maineqnupper} \E\lmax(\Sigma_N) \le1+\varepsilon.
\end{equation}
Here $C_\upper:= 512 (16C)^{1+2/\eta} (6+6/\eta)^{1+4/\eta}$.

\end{theorem}

We shall control the largest eigenvalue of a symmetric matrix $A$
using the (upper) \textit{Stieltjes transform}
\[
\mp_A(u) = \tr(uI-A)^{-1} = \sum
_{i=1}^n \bigl( u-\lambda_i(A)
\bigr)^{-1},\qquad u \in\R.
\]
Similarly to our argument for the lower edge, for a sensitivity value
$\pn> 0$,
we define the \textit{upper soft spectral edge} $u_\pn(A)$ to be the
largest $u$ for which
\[
\mp_A(u) = \pn.
\]
Since $\mp_A(u)$ decreases from $\infty$ to $0$ as $u$ increases from
the upper
spectral edge $\lmax(A)$ to $\infty$, the value $u_\pn(A)$ is
defined uniquely, and
\[
u_\pn(A) > \lmax(A).
\]
For $\pn\to\infty$ we have $u_\pn(A) \to\lmax(A)$, but as before,
we shall work with small sensitivity values $\pn\in(0,1)$.
Our goal is to show that $u_\pn(A)$ increases by about $1$, on
average, with every rank-one update.

%th3.2 #&#
\begin{theorem}[(Random upper shift)] \label{upperrankone}
Suppose $X$ is an isotropic random vector satisfying the strong regularity
assumption (\SR) for some $C,\eta> 0$.
Assume $\varepsilon\in(0,1)$ and
%
%e3.2 #&#
\begin{equation}
\label{upperbarriergamma}
\pn\le c_{\rref{upperrankone}} \varepsilon^{1+2/\eta},
\end{equation}
where $c_{\rref{upperrankone}}^{-1} = 256 (8C)^{1+2/\eta} (6+6/\eta
)^{1+4/\eta}$.
Then for every symmetric matrix~$A$, one has
%
%e3.3 #&#
\begin{equation}
\label{eqnmainupper} \E u_\pn\bigl(A+XX^T\bigr) \le
u_\pn(A) + 1+\varepsilon.
\end{equation}
\end{theorem}

Iterating Theorem~\ref{upperrankone} yields a proof of Theorem \ref
{mainupper}.
\begin{pf*}{Proof of Theorem~\ref{mainupper}}
The argument is similar to the proof of Theorem~\ref{mainlower} given
in Section~\ref{sloweredge}.
We set $\pn= \pn(\varepsilon) = c_{\rref{upperrankone}}
\varepsilon^{1+2/\eta}$.
Then we start with $A_0=0$ where $u_\pn(A_0) = n/\pn$, and we
inductively apply Theorem~\ref{upperrankone}
for $A_k=A_{k-1}+X_kX_k^T$ to obtain
\[
\E\lmax \Biggl( \frac{1}{N}\sum_{i=1}^N
X_iX_i^T \Biggr) < \frac{u_\pn(A_0)}{N} + 1+
\varepsilon= 1 + \varepsilon+ \frac
{n}{\pn N}.
\]
For $N \ge n/\varepsilon\pn$, the bound becomes $1+2\varepsilon$.
Substituting the value of $\pn$ and replacing $\varepsilon$ by
$\varepsilon/2$ gives the promised result.\vadjust{\goodbreak}
\end{pf*}

The above proof works for $\varepsilon,\pn(\varepsilon) < 1$ and
thus for $N = \Omega(n)$,
but it may be extended to smaller $N$ as follows.

\begin{pf*}{Proof of Corollary~\ref{fixedN}}
In the proof of Theorem~\ref{mainupper}, we have shown that for
every $\varepsilon\in(0,1)$ and every positive integer $N$, we have
\[
E:= \E\lmax(\Sigma_N) < 1 + \varepsilon+ \frac{n}{\pn
(\varepsilon) N},
\]
where $\pn(\varepsilon) = c_{\rref{upperrankone}}  \varepsilon^{1+2/\eta}$. Optimizing in $\varepsilon$, we
apply this estimate with $\varepsilon= (n/N)^{{1}/({2+2/\eta})}$
when $n<N$ and with
$\varepsilon=1/2$ when $n \ge N$ to obtain
\begin{eqnarray*}
E &<& 1 + \bigl(1+c_{\rref{upperrankone}}^{-1}\bigr) \biggl(\frac{n}{N}
\biggr)^{{1}/({2+2/\eta})} \qquad\mbox{if } n<N,
\\
E &<& \frac{3}{2} + \frac{n}{\pn(1/2)N} \le1 + 2^{2+2/\eta}c_{\rref{upperrankone}}^{-1}
\biggl(\frac{n}{N} \biggr) \qquad\mbox{if } n \ge N.
\end{eqnarray*}
Combining these, for every $n$ and $N$ we conclude that
\[
E < 1 + \bigl(1+c_{\rref{upperrankone}}^{-1}\bigr) \biggl(\frac{n}{N}
\biggr)^{
{1}/({2+2/\eta})} + 2^{2+2/\eta}c_{\rref{upperrankone}}^{-1} \biggl(
\frac{n}{N} \biggr)
\]
as required.

A similar bound for $\E\lmin(\Sigma_N)$ is immediate from Theorem
\ref{mainlower}; see the remark after its proof.
\end{pf*}

%that the proofs of Theorem~\ref{mainupper}
%and Corollary~\ref{fixedN} do not actually require the independence of
%$X_i$.
%To invoke Theorem~\ref{upperrankone}, we only need that the
%distribution
%of each $X_i$ conditioned on $X_1,\ldots, X_{i-1}$ satisfies

The rest of this section is devoted to proving Theorem~\ref{upperrankone}.
Given a matrix~$A$, a real number $u>\lmax(A)$ and a vector $x \in\R^n$,
we say that $\D\ge0$ is a \textit{feasible upper shift} if
%
%e3.4 #&#
\begin{equation}
\label{upperfeasibility} A + xx^T \prec(u+\D)I
\quad\mbox{and}\quad
\mp_{A+xx^T}(u+\D) \le\mp_A(u).
\end{equation}
The definition of the soft spectral edge $u = u_\pn(A)$ along with
monotonicity of the Stieltjes transform implies
that
%
%e3.5 #&#
\begin{equation}
\label{spectraledgeuppershift} u_\pn\bigl(A+xx^T
\bigr) \le u_\pn(A) + \D
\end{equation}
for every feasible upper shift $\D$. So will be done if we can produce
a feasible shift $\D$
such that $\E\D\le1+\varepsilon$ where the expectation is over
random $X$.

As in our argument for the lower edge, we begin by reducing the
feasibility for a shift $\dn$ to an
inequality involving two quadratic forms.

%le3.3 #&#
\begin{lemma}[(Feasible upper shift)] \label{upperfeasible}
Consider the numbers $u \in\R$, $\D> 0$, a~matrix $A \prec uI$ and a
vector $x$.
Then a sufficient condition for $\D\ge0$ to be a feasible upper shift is
%
%e3.6 #&#
\begin{eqnarray}
\label{eupperfeasible}
&&\frac{x^T(u+\D-A)^{-2}x}{\mp_A(u)-\mp_A(u+\D)} + x^T(u+
\D-A)^{-1}x \nonumber\\[-8pt]\\[-8pt]
&&\qquad=: \Qp_2(\D,x) + \Qp_1(\D,x) \le1.\nonumber
\end{eqnarray}
\end{lemma}

\begin{pf}
Note that $A\prec uI \prec(u+\D) I$ so that all quadratic forms are positive,
and assume $x\ne0$ since otherwise the claim is trivial.
As in the proof of Lem\-ma~\ref{lowerfeasible}, we use the Sherman--Morisson
formula to write
\begin{eqnarray*}
\mp_{A+xx^T}(u+\D) &=& \tr\bigl(u+\D-A-xx^T
\bigr)^{-1}
\\
&=& \mp_{A}(u+\D) + \frac{x^T(u+\D-A)^{-2}x}{1-x^T(u+\D-A)^{-1}x}
\\
&=& \mp_{A}(u) - \bigl(\mp_A(u)-\mp_A(u+\D)
\bigr)\\
&&{} + \frac{x^T(u+\D-A)^{-2}x}{1-x^T(u+\D-A)^{-1}x}.
\end{eqnarray*}
Rearranging reveals that $\mp_{A+xx^T}(u+\D)\le\mp_A(u)$ exactly when
(\ref{eupperfeasible}) holds.

To establish the second condition
%
%e3.7 #&#
\begin{equation}
\label{noupperjump} xx^T \prec u+\D-A,
\end{equation}
we recall that
\[
R\prec S \quad\iff\quad S^{-1/2}RS^{-1/2} \prec I
\]
for all positive matrices $R,S$ (this can be seen, e.g., using the
Courant--Fischer theorem). Applying this fact to (\ref{noupperjump}),
we see that
it suffices to have
\[
(u+\D-A)^{-1/2}xx^T (u+\D-A)^{-1/2} \prec I
\]
or equivalently
\[
x^T (u+\D-A)^{-1} x < 1,
\]
which follows from (\ref{eupperfeasible}) and $\Qp_2(\D,x)>0$.
\end{pf}

We will reason about the two quantities $\Qp_1$ and $\Qp_2$ separately,
producing two separate shifts $\D_1$ and $\D_2$ for them and eventually
combining these into a single $\D:=\D_1 \lor\D_2$, as required by
Lemma~\ref{upperfeasible}.

For some fixed parameter $\tau\in(0,1)$, let us define $\D_1=\D_1(A,x,u)$ and $\D_2=\D_2(A,x,u)$
to be the smallest nonnegative numbers such which satisfy
%
%e3.8 #&#
\begin{equation}
\label{Q1Q2bounds} \Qp_1(\D_1,x)\le\tau,\qquad
\Qp_2(\D_2,x)\le1-\tau.
\end{equation}
For $u=u_\pn(A)$ and for a random vector $x=X$,
Lemmas~\ref{controlq1} and~\ref{controlq2} will allow us to control
the expected
value of each of these shifts, so
%
%e3.9 #&#
\begin{equation}
\label{expecteduppershifts} \E\D_1 \le\varepsilon/2, \qquad\E
\D_2 \le1+\varepsilon/2,
\end{equation}
whenever the sensitivity parameter $\pn=\pn(\tau,\varepsilon)$ is
sufficiently
small.
From this we will obtain Theorem~\ref{upperrankone} quickly as follows.

\begin{pf*}{Proof of Theorem~\ref{upperrankone}}
Let $u_\pn(A)=u$, so the condition $A \prec uI$ of Lem\-ma~\ref{upperfeasible} holds.
Consider the shifts $\D_1=\D_1(A,X,u)$ and $\D_2=\D_2(A,X,u)$
defined above.
By (\ref{Q1Q2bounds}), we have
\[
Q_1(\D_1,X) + Q_2(\D_2,X) \le1.
\]
Moreover, a quick inspection of the quadratic forms in Lemma
\ref{upperfeasible} shows
that $Q_1(\D,X)$ and $Q_2(\D,X)$ are decreasing in $\D$, and hence
\[
Q_1(\D_1 \vee\D_2,X) + Q_2(
\D_1 \vee\D_2,X) \le1.
\]
Then Lemma~\ref{upperfeasible} guarantees that $\D_1 \vee\D_2$ is a
feasible upper shift,
which implies by (\ref{spectraledgeuppershift}) that
\[
u_\pn\bigl(A+XX^T\bigr) \le u_\pn(A) +
\D_1 \vee\D_2.
\]
Furthermore, (\ref{expecteduppershifts}) yields a bound on the
expected shift
\[
\E\D_1\lor\D_2 \le\E\D_1+\E\D_2
\le1+\varepsilon,
\]
which gives conclusion (\ref{eqnmainupper}) of Theorem~\ref{upperrankone}.

It remains to note that Lemmas~\ref{controlq1} and~\ref{controlq2}
only guarantee
that the bounds (\ref{expecteduppershifts}) hold when the sensitivity
$\pn$
is sufficiently small, namely $\pn\le\pn_1(\tau,\varepsilon/2)
\wedge\pn_2(\tau,\varepsilon/2)$.
With $\tau= \varepsilon/16$, we can simplify this inequality into
the assumption of
Theorem~\ref{upperrankone}.
\end{pf*}

The rest of this section is devoted to controlling the shifts $\D_1$
and $\D_2$.

\begin{remark*} It is easy to check that the proofs of Lemmas \ref
{controlq1} and~\ref{controlq2} which follow, and consequently Theorem~\ref{upperrankone}, only require
%
%e3.10 #&#
\begin{equation}
\label{weakisotropy2} \E X_iX_i^T\prec c I
\end{equation}
for some constant $c\!=\!c(\varepsilon)$.
Thus
if we desire a bound of $\lmax (\frac{1}{N}\sum_{i=1}^N
X_iX_i^T )\!<
1+\varepsilon$ in Theorem~\ref{mainupper}, then $\E X_iX_i^T = I$
can be
replaced by the weaker condition (\ref{weakisotropy2}).
\end{remark*}

%s3.1 #&#
\subsection{\texorpdfstring{Control of $\D_1$}{Control of Delta 1}}

%le3.4 #&#
\begin{lemma} \label{controlq1}
Consider numbers $u \in\R$, $\pn> 0$ and a matrix $A \prec uI$
satisfying $\mp_A(u) \le\pn$.
Let $X$ be a random vector satisfying (\ref{allmarginals}) for some
$C,\eta>0$,
and let $\varepsilon, \tau\in(0,1)$.
If the sensitivity satisfies
\[
\pn\le\pn_1(\tau,\varepsilon):= \frac{\tau^{1+1/\eta} \varepsilon^{1/\eta}}{(4C)^{1+1/\eta
}(4+4/\eta)^{1+3/\eta}},
\]
then the shift $\D_1 = \D_1(A,X,u)$ satisfies
\[
\E\D_1 \le\varepsilon.
\]
\end{lemma}

\begin{pf}
Let $(\psi_i)_{i\le n}$ and $(\l_i)_{i\le n}$ denote the eigenvectors
and eigenvalues of~$A$,
and let $\xi_i = \langle X,\psi_i\rangle^2$.
We know that $\mp_A(u) = \sum_{i=1}^n (u-\lambda_i)^{-1} \le\pn$,
and $\D_1$ is the smallest nonnegative number satisfying
\[
\sum_{i=1}^n \frac{\xi_i}{u-\lambda_i+\D_1}\le\tau.
\]
Rescaling everything by $\pn$ and setting $\mu_i:= \pn(u-\lambda_i)$ so that
\[
\sum_{i=1}^n \frac{1}{\mu_i} = \sum
_{i=1}^n \frac{1}{\pn
(u-\lambda_i)} \le1,
\]
the problem becomes equivalent to bounding the least
$\mu:= \pn\D_1$ for which
\[
\sum_{i=1}^n \frac{1}{\mu_i+\mu}\le
\frac{\tau}{\pn}.
\]
Applying the following, somewhat more general, probabilistic lemma to
$(\xi_i)_{i\le n}$, we conclude that
\[
\E\D_1 \le\frac{1}{\pn}\E\mu\le \frac{1}{\pn}
\frac{C(4+4/\eta)^{3+\eta}(4\pn)^{1+\eta}}{\tau^{1+\eta}},
\]
whenever
\[
\pn\le\frac{\tau}{4C}.
\]
Substituting $\pn=\pn_1(\tau,\varepsilon)$ gives the promised bound.
\end{pf}

%le3.5 #&#
\begin{lemma}
Suppose $\{\xi_i\}_{i\le n}$ are positive random variables with $\E
\xi_i = 1$
and
%
%e3.11 #&#
\begin{equation}
\label{allmarginalssubsets} \P \biggl\{ \sum_{i\in S}
\xi_i \ge t \biggr\}\le \frac{C}{t^{1+\eta}}
\quad\mbox{provided}\quad t>
C|S|=C\sum_{i\in S} \E\xi_i
\end{equation}
for all subsets $S\subset[n]$ and some constants $C,\eta>0$.
Consider positive numbers $\mu_i$ such that
\[
\sum_{i=1}^n \frac{1}{\mu_i} \le1.
\]
Let $\mu$ be the minimal positive number such that
\[
\sum_{i=1}^n \frac{\xi_i}{\mu_i+\mu} \le K
\]
for some $K\ge4C$.
Then
\[
\E\mu\le\frac{C(4+4/\eta)^{3+\eta}}{(K/4)^{1+\eta}}.
\]
\end{lemma}

\begin{pf}
For simplicity of calculations, assume for the moment that the values
of all $\mu_i$
are dyadic, that is,
\[
\mu_i \in\bigl\{ 2^0, 2^1, 2^2,
\ldots\bigr\}.
\]
For each dyadic number $k$, let
\[
I_k:= \{ i\dvtx u_i=k \},\qquad n_k:=
|I_k|.
\]
By assumption, we have
\[
1 \ge\sum_{i=1}^n \frac{1}{\mu_i} =
\sum_{k\ \mathrm{dyadic}} \sum_{i \in I_k}
\frac{1}{k} = \sum_{k\ \mathrm{dyadic}} \frac{n_k}{k},
\]
and $\mu$ is the smallest positive number such that
%
%e3.12 #&#
\begin{equation}
\label{mu} \sum_{i=1}^n
\frac{\xi_i}{\mu_i+\mu} = \sum_{k\ \mathrm{dyadic}} \frac{1}{k+\mu}
\sum_{i \in I_k} \xi_i \le K.
\end{equation}

We estimate $\mu$ by replacing it with a bigger but easier quantity
$\mu'$.
Define $\mu'$ to be the smallest positive number such that, for every
dyadic $k$, one has
\[
\frac{1}{k+\mu'} \sum_{i \in I_k} \xi_i \le
\varepsilon_k \qquad\mbox{where } \varepsilon_k:=
\frac{K}{2} \frac{n_k}{k} \lor \frac{K}{2\sigma}
k^{-{\eta}/({2+2\eta})},
\]
where
%
%e3.13 #&#
\begin{equation}
\label{sigma} \sigma:= \sum_{\mathrm{dyadic}\ k} k^{-{\eta}/({2+2\eta})}
\le\frac{2+2\eta}{\eta}\sum_{\mathrm{dyadic}\ k} \frac{1}{k}
\le4+4/\eta.
\end{equation}
Since
\[
\sum_{k\ \mathrm{dyadic}} \frac{1}{k+\mu'} \sum
_{i \in I_k} \xi_i \le\sum_{k\ \mathrm{dyadic}}
\varepsilon_k \le\frac{K}{2} \sum
_{k\ \mathrm{dyadic}} \frac{n_k}{k} + \frac
{K}{2\sigma} \sum
_{k\ \mathrm{dyadic}} k^{-{\eta}/({2+2\eta})} \le K,
\]
the definition of $\mu$ given in (\ref{mu}) yields
\[
\mu\le\mu'.
\]
It remains to bound $ \E\mu'$.

By definition,
\[
\mu' = \max_{k\ \mathrm{dyadic}} \biggl( \frac{1}{\varepsilon_k} \sum
_{i \in I_k} \xi_i - k \biggr)_+.
\]
Let $\theta_k = \frac{1}{\varepsilon_k} \sum_{i \in I_k} \xi_i -
k$. For every $t \ge0$, one has
\[
\P\{ \theta_k > t \} = \P \biggl\{ \sum
_{i \in I_k} \xi_i > (k+t) \varepsilon_k
\biggr\}.
\]
Since $\varepsilon_k \ge\frac{K n_k}{2k}$ by definition, we have
\[
(k+t)\varepsilon_k \ge k \varepsilon_k \ge
\frac{K n_k}{2} = \frac
{K}{2} \E \biggl( \sum
_{i \in I_k} \xi_i \biggr) \ge C \E \biggl( \sum
_{i \in I_k} \xi_i \biggr).
\]
So by regularity assumption (\ref{allmarginalssubsets}),
\[
\P\{ \theta_k > t \} \le\frac{C}{(k+t)^{1+\eta} \varepsilon_k^{1+\eta}}.
\]
A union bound then gives
\begin{eqnarray*}
\P\bigl\{ \mu' > t \bigr\} &\le& \sum_{k\ \mathrm{dyadic}}
\frac{C}{(k+t)^{1+\eta} \varepsilon_k^{1+\eta}}
\\
&\le& \frac{C}{(K/2\sigma)^{1+\eta}} \sum_{k\ \mathrm{dyadic}} \frac
{k^{\eta/2}}{(k+t)^{1+\eta}}
\qquad\mbox{(by definition of $\varepsilon_k$)}
\\
&\le& \frac{C}{(K/2\sigma)^{1+\eta}} \sum_{k\ \mathrm{dyadic}}
\frac
{1}{(k+t)^{1+\eta/2}}.
\end{eqnarray*}
This implies that
\begin{eqnarray*}
\E\mu' &=& \int_0^\infty\P \bigl\{
\mu'>t \bigr\} \,dt \le\frac{C}{(K/2\sigma)^{1+\eta}}\sum
_{k\ \mathrm{dyadic}} \int_0^\infty
\frac{dt}{(k+t)^{1+\eta/2}}
\\
&=& \frac{C}{(K/2\sigma)^{1+\eta}} \sum_{k\ \mathrm{dyadic}} \frac
{k^{-\eta/2}}{\eta/2}
\\
&\le& \frac{C}{(K/2\sigma)^{1+\eta}} \frac{2}{\eta}\cdot\frac{4}{\eta}
\qquad\mbox{[by a calculation similar to (\ref{sigma})]}
\\
&\le& \frac{C}{(K/2)^{1+\eta}}(4+4/\eta)^{3+\eta} \qquad\mbox{[by (\ref{sigma})]}.
\end{eqnarray*}
The promised bound for general (nondyadic) $\mu_i$ follows by
rounding each
$\mu_i$ down to the nearest power of $2$ and replacing $K$ by $K/2$.
\end{pf}

\begin{remark*}[{[Necessity of the strong regularity assumption (\ref
{allmarginals})]}]
The preceding lemma is the only place in the proof where the full power
of (\ref{allmarginals}) is used.
To see that it is necessary, consider the following
situation. Fix any $S\subset[n]$, and let
$\frac1{\mu_i} = \one_{\{i\in S\}}|S|$ so that $\sum_i \frac
{1}{\mu_i}=1$. Then the\vspace*{2pt}
smallest $\mu\ge0$ for which $\sum_i \frac{1}{\mu_i+\mu}\le K$ is just
\[
\mu= \biggl(\frac{1}{K} \sum_{i\in S}
\xi_i - |S| \biggr)_+.
\]
We now lowerbound the tail probability
\[
\P\{\mu\ge t\} = \P \biggl\{\sum_{i\in S}
\xi_i \ge K\bigl(|S|+t\bigr) \biggr\} \ge\P \biggl\{\sum
_{i\in S}\xi_i \ge2Kt \biggr\} \qquad\mbox{for } t \ge|S|.
\]
In order to have $\E\mu= O(1)$, this probability must be
$O(1/t)$ by Markov's inequality, which is essentially assumption (\ref
{allmarginalssubsets}) of the lemma.
In the proof of Theorems~\ref{main} and~\ref{upperrankone},
the sums of random variables $\xi_i$ arise from projections of the
random vector $X$ onto varying
eigenspaces of $A$; the only succinct way to guarantee (\ref{allmarginalssubsets}) for all such
projections is essentially (\ref{allmarginals}).
\end{remark*}

%s3.2 #&#
\subsection{\texorpdfstring{Control of $\D_2$}{Control of Delta 2}}

%le3.6 #&#
\begin{lemma} \label{controlq2}
Consider numbers $u \in\R$, $\pn> 0$ and a matrix $A \prec uI$
satisfying $\mp_A(u) \le\pn$.
Let $X$ be a random vector satisfying (\ref{allmarginals}) for some
$C,\eta>0$,
and let $\varepsilon\in(0,1)$, $0 < \tau< \varepsilon/2$ be parameters.
If the sensitivity satisfies
\[
\pn\le\pn_2(\tau,\varepsilon):= \frac{\varepsilon^{2/\eta}(\varepsilon-4\tau)}{128\cdot
(2C)^{2/\eta} (4+6/\eta)^{4/\eta}},
\]
then the shift $\D_2 = \D_2(A,X,u)$ satisfies
\[
\E\D_2 \le1 + \varepsilon.
\]
\end{lemma}

It will be more convenient to work with the quadratic form
\[
\Qp_2'(\D,x):= \frac{x^T(u+\D-A)^{-2}x}{\tr(u+\D-A)^{-2}},
\]
for which we have
%
%e3.14 #&#
\begin{equation}
\label{Q2Q2} \frac{1}{\D}\Qp_2'(
\D,x)\ge\Qp_2(\D,x) \qquad\mbox{for } \D> 0,
\end{equation}
since the denominators satisfy
\[
\mp_A(u)-\mp_A(u+\D) = \tr\bigl[(uI-A)^{-1}
- (u+\D-A)^{-1}\bigr] \ge\D\tr(u+\D-A)^{-2}.
\]

\begin{remark*} The reason for working with $\Qp_2$ rather than
directly with $\Qp_2'$
in Lemma~\ref{upperfeasible} is that $\Qp_2(\D,x)$ is decreasing in
$\D$; this monotonicity is required when
arguing that the maximum of the two shifts $\D=\D_1\lor\D_2$ is
feasible in the proof of Theorem~\ref{upperrankone}.
\end{remark*}

We begin by recording some regularity properties of $\Qp_2'(\D,X)$.

%le3.7 #&#
\begin{lemma}[{[Regularity and moments of of $\Qp_2'(\D,X)$]}] \label
{regularityf}
Consider numbers $u \in\R$, $\pn> 0$ and a matrix $A \prec uI$
satisfying $\mp_A(u) \le\pn$.
Let $X$ be a random vector satisfying (\ref{allmarginals}) for some
$C,\eta>0$.
Then for every $\D\ge0$ one has:
\begin{longlist}
\item$\Qp_2'(\D,X) \le(1+\pn\D)^2  \Qp_2'(0,X)$;
\item$\E\Qp_2'(\D,X) = 1$;
\item$\E\Qp_2'(\D,X)^p \le C(3+3/\eta)$ for $p=1+2\eta/3$.
\end{longlist}
\end{lemma}

\begin{pf}
(i) is analogous to Lemma~\ref{lemlowerquadraticregularity}. In a
similar way, we show that
all eigenvalues $\l_i$ of $A$ satisfy $u-\l_i \ge1/\pn$, which
implies the comparison inequality
\[
u-\l_i \le u+\D- \l_i \le(1+\pn\D) (u-
\l_i).
\]
Denoting $(\psi_i)_{i\le n}$ the eigenvectors of $A$, we express
%
%e3.15 #&#
\begin{equation}
\label{Q2sum} Q_2'(\D,X) =
\frac{\sum_{i=1}^n (u+\D-\lambda_i)^{-2} \langle
X,\psi_i\rangle^2} {
\sum_{i=1}^n (u+\D-\lambda_i)^{-2}}.
\end{equation}
The comparison inequality yields (i).

\mbox{}\hphantom{i}(ii) We note that (\ref{Q2sum}) can be rearranged as a convex
combination of $\langle X,\psi_i\rangle^2$.
\[
Q_2'(\D,X) = \sum_i
\alpha_i \langle X,\psi_i\rangle^2
\qquad\mbox{where } \alpha_i \ge0, \sum_{i=n}
\alpha_i = 1.
\]
Then (ii) follows since $\E\langle X,\psi_i\rangle^2 = 1$ by isotropy.

(iii) We apply Minkowski's inequality to obtain
\[
\bigl(\E Q_2'(\D,X)^p\bigr)^{1/p}
\le\sum_{i=1}^n \alpha_i
\bigl(\E\langle X,\psi_i\rangle^{2p}\bigr)^{1/p}.
\]
Now a simple integration of tails implies that each
\[
\E\langle X,\psi_i\rangle^{2p} = \E\langle X,
\psi_i\rangle^{2+4\eta/3} \le C(3+3/\eta),
\]
which concludes the proof.
\end{pf}
Next, we see how the regularity properties of $\Qp_2'(\D,X)$
translate into the corresponding
properties of $\D_2$:

%le3.8 #&#
\begin{lemma}[(Regularity of $\D_2$)] \label{regularitym}
Consider numbers $u \in\R$, $\pn> 0$ and a matrix $A \prec uI$
satisfying $\mp_A(u) \le\pn$.
Let $X$ be a random vector satisfying (\ref{allmarginals}) for some
$C,\eta>0$,
and let $0 < \tau< 1/2$.
Then the shift $\D_2 = \D_2(A,X,u)$ satisfies:
\begin{longlist}
\item$\E\D_2^{1+\eta/2} \le2^{1+\eta}C(4+6/\eta)^2$;
\item$\E\D_2 \one_{\{\Qp_2'(0,X) \le(t-2\tau)/8\pn\}} \le1 +
t$ for every $t \in[0,1]$.
\end{longlist}
\end{lemma}

\begin{pf}
(i) By definition of $\D_2$ and using (\ref{Q2Q2}), we have for all $t>0$,
\[
\P\{\D_2>t\} \le\P\bigl\{\Qp_2(t,X) > 1-\tau\bigr\} \le
\P\bigl\{\Qp_2'(t,X)>t(1-\tau)\bigr\}.\vadjust{\goodbreak}
\]
This probability can be controlled using
Lemma~\ref{regularityf}(iii) and Markov's inequality,
so we obtain
\[
\P\{\D_2>t\} \le\frac{C(3+3/\eta)}{t^{1+2\eta/3}(1-\tau
)^{1+2\eta/3}} \le\frac{C(3+3/\eta)}{(1/2)^{1+2\eta/3}t^{1+2\eta/3}}
\]
as $\tau<1/2$. Integration of tails yields
\[
\E\D_2^{1+\eta/2} \le2^{1+2\eta/3}\cdot C(3+3/\eta) (4+6/\eta
),
\]
which implies the claim.

(ii) Let $s_0$ denote the smaller solution of the quadratic equation
\[
(1+s\pn)^2 \Qp_2'(0,X) = s(1-\tau),
\]
whenever a solution exists. In this case $s_0>0$ and
Lemma~\ref{regularityf}(i) yields that
\[
\Qp_2'(s_0,X) \le s_0(1-\tau).
\]
By (\ref{Q2Q2}), this yields $\Qp_2(s_0,X) \le s_0(1-\tau)$.
By definition of $\D_2$, this in turn implies that
\[
\D_2 \le s_0.
\]
An elementary calculation shows that if $\Qp_2'(0,X) \le(t-2\tau
)/8\pn$, then
the solution $s_0$ exists and satisfies
\[
s_0 \le(1+t ) \Qp_2'(0,X).
\]
It follows that
\[
\E s_0 \one_{\{\Qp_2'(0,X) \le(t-2\tau) /8\pn\}} \le(1+t ) \E\Qp_2'(0,X)
= 1+t,
\]
where we used Lemma~\ref{regularityf}(i) in the last step.
\end{pf}

We can now complete the proof of Lemma~\ref{controlq2}.

\begin{pf*}{Proof of Lemma~\ref{controlq2}}
We decompose
\[
\E\D_2 = \E\D_2 \one_{\{\Qp_2'(0,X) \le(t-2\tau)/8\pn\}} + \E
\D_2 \one_{\{\Qp_2'(0,X) > (t-2\tau)/8\pn\}} =: E_1 + E_2.
\]
By Lemma~\ref{regularitym}(ii), we have $E_1 \le1+t$.
Next, we estimate $E_2$ using H\"older's inequality,
\[
E_2 \le \bigl( \E\D_2^{1+\eta/2}
\bigr)^{{1}/({1+\eta/2})} \bigl( \P\bigl\{\Qp_2'(0,X)>(t-2
\tau)/8\pn\bigr\} \bigr)^{({\eta
/2})/({1+\eta/2})}.
\]
The two terms here can be estimated using Lemma~\ref{regularitym}(i)
and Lemma~\ref{regularityf} along with Markov's inequality,
\begin{eqnarray*}
E_2 &\le& \bigl(2^{1+\eta}C(4+6/\eta)^2
\bigr)^{{1}/({1+\eta/2})} \biggl(\frac{C(3+3/\eta)}{((t-2\tau)/8\pn)^{1+\eta/2}} \biggr)
^{({\eta/2})/({1+\eta/2})}
\\
&\le& 2^{1+\eta}C(4+6/\eta)^2\cdot \biggl(\frac{8\pn}{t-2\tau}
\biggr)^{\eta/2}.
\end{eqnarray*}
Finally, we set $t=\varepsilon/2$ and use the assumptions $\pn\le
\pn_2(\tau,\varepsilon)$ and
$\t< \varepsilon/2$ to conclude that $E_2 \le\varepsilon/2$.
Together with $E_1 \le1+ t = 1+\varepsilon/2$
this implies
\[
\E\D_2 \le1+\varepsilon
\]
as claimed.
\end{pf*}

\begin{remark*}
Although for convenience of application Lemma~\ref{controlq2} is stated
under the strong regularity assumption (\ref{allmarginals}),
the latter is not used in the proof. The argument above uses only
the weak regularity assumption (\ref{onedimmarginals}).
\end{remark*}

%The choice of $\pn$ is dictated by a trade-off between the total
%number of steps
%we plan to take $N$, and the precision $\e$ with which we want to
%control
%the spectrum. Smaller $\pn$ will imply that $\ell_\pn(A_k)$ is more
%stable and
%that we can take the expected shifts $\E(\ell_\pn(A_k)-\ell_
%be at least $1-\e$ for smaller $\e$, eventually leading to
%closer control of $\lmin(A_N)$; at the same time it will
%force us to take the initial barrier $\ell_0=-n/\pn$ to be very
%negative, which
%will require larger
%number of steps $N$ to overwhelm into the same magnitude as the error
%term
%$\e$ (roughly $n/\e\pn$ steps).

%s4 #&#
\section{The spectral norm} \label{sspectralnorm}
In this section we prove Theorem~\ref{main} by showing that whenever
$X_1,\ldots,X_N$ are independent and satisfy (\ref{allmarginals}),
the spectral
norm estimate
%
%e4.1 #&#
\begin{equation}
\label{normbound} \E\| \Sigma_N - I \| \le\varepsilon
\end{equation}
follows from the spectral edge estimates
%
%e4.2 #&#
\begin{equation}
\label{edgebounds} \E\lmin(\Sigma_N) \ge1-\varepsilon/3;\qquad \E\lmax (
\Sigma_N) \le1+\varepsilon/3
\end{equation}
obtained in Theorems~\ref{mainlower} and~\ref{mainupper}.
The basic idea is to show using independence that
\[
\lambda_{\mathrm{average}}(\Sigma_N) = \frac{1}{n}\tr(
\Sigma_N)
\]
is concentrated near its expectation of $1$.
Combining this with
\[
\E \bigl(\lmax(\Sigma_N) - \lmin(\Sigma_N) \bigr) \le2
\varepsilon /3,
\]
which follows immediately from (\ref{edgebounds}), yields (\ref{normbound}).

We rely on the following elementary proposition regarding
sums of independent random variables.

%pr4.1 #&#
\begin{proposition}\label{indepprop}
Let $Z_i$ be independent random variables with $\E Z_i = 1$ and
satisfying the following tail bounds for some $C,\eta>0$:
\[
\P\bigl\{ |Z_i| > t \bigr\} \le C t^{-1-\eta}, \qquad t > 0.
\]
If $\varepsilon\in(0,1)$ and
\[
N \ge\frac{(2C)^{2/\eta}(1+1/\eta)^{2/\eta}}{(\varepsilon
/2)^{2+2/\eta}},
\]
then
\[
\E \Biggl| \frac{1}{N} \sum_{i=1}^N
Z_i -1 \Biggr| \le\varepsilon.
\]
\end{proposition}

Postponing the proof of Proposition~\ref{indepprop}, we use this fact
to control
\[
\frac{1}{n}\tr(\Sigma_N) = \frac{1}{n} \sum
_{i=1}^N \frac{\|X_i\|_2^2}{N}
\]
and prove the main theorem as follows.

\begin{pf*}{Proof of Theorem~\ref{main}}
Assume the random vectors $X_i$ are isotropic and satisfy (\ref{allmarginals})
with parameters $C,\eta$. This implies that the random variables
\[
Z_i = \frac{\|X_i\|_2^2}{n}
\]
satisfy the requirements of Proposition~\ref{indepprop} with parameters
$C^{1+\eta},\eta$.
It follows that
%
%e4.3 #&#
\begin{equation}
\label{tracebound}
\E \biggl|\frac{1}{n} \tr(\Sigma_N - I) \biggr| = \E
\Biggl| \frac {1}{N} \sum_{i=1}^N Z_i -1 \Biggr| \le\varepsilon,
\end{equation}
whenever
%
%e4.4 #&#
\begin{equation}
\label{Nbig} N\ge\frac{(4C)^{2+2/\eta}(1+1/\eta
)^{2/\eta}}{\varepsilon^{2+2/\eta}} =: \frac{C_\prop}{\varepsilon^{2+2/\eta}}.
\end{equation}

Now consider the random variables
\[
L = \lmin(\Sigma_N - I),\qquad U = \lmax(\Sigma_N - I),\qquad M =
\frac{1}{n} \tr(\Sigma_N - I).
\]
We have
\[
L \le M \le U,
\]
and we are interested in
%
%e4.5 #&#
\begin{equation}
\label{ulm} \|\Sigma_N - I\| = U \vee-L \le U - L + |M|.
\end{equation}

When $N\ge C_\upper n/\varepsilon^{2+2/\eta}$, Theorem~\ref{mainupper} gives $\E U \le\varepsilon$.
To show that \mbox{$\E L \ge\varepsilon$},
we recall that (\ref{allmarginals}) with parameters $C,\eta$ implies
(\ref{onedimmarginals}) with parameters $C(2+2/\eta),\eta$ and
invoke Theorem
\ref{mainlower},
noting that its requirement (\ref{Cprime2}) is satisfied as
\[
C_\upper=512 (16C)^{1+2/\eta} (6+6/\eta)^{1+4/\eta} > 40
\bigl(10C(2+2/\eta)\bigr)^{2/\eta} = C_\low.
\]
Now that we have both bounds $\E U \le\varepsilon$ and $\E L \ge
\varepsilon$, we can combine them
with (\ref{tracebound}) and (\ref{ulm}), which yields
\[
\E\|\Sigma_N - I\| \le2\varepsilon+ \varepsilon,
\]
whenever
%
%e4.6 #&#
\begin{equation}
\label{finalopt} N \ge C_\upper\frac{n}{\varepsilon^{2+2/\eta}} \lor
C_\prop\frac
{1}{\varepsilon^{2+2/\eta}}.
\end{equation}
Replacing $\varepsilon$ by $\varepsilon/3$ and taking
\[
N\ge C_\main\frac{n}{\varepsilon^{2+2/\eta}},
\]
where
\[
C_\main:= 512 \cdot3^{2+2/\eta}\cdot(16C)^{2+2/\eta} (6+6/\eta
)^{1+4/\eta}
\]
always satisfies (\ref{finalopt}). This completes the proof of the theorem.
\end{pf*}

\begin{pf*}{Proof of Proposition~\ref{indepprop}}
Fix a parameter $K>0$, and decompose
\[
Z_i = Z_i \one_{\{|Z_i| \le K\}} + Z_i
\one_{\{|Z_i| > K\}} =: Z_i' + Z_i''.
\]
Using $\E Z_i' + \E Z_i'' = \E Z_i = 1$ and by triangle inequality, we obtain
\begin{eqnarray*}
\E \Biggl| \frac{1}{N} \sum_{i=1}^N
Z_i -1 \Biggr| &\le&\E \Biggl| \frac{1}{N} \sum
_{i=1}^N Z_i' - \E
\frac{1}{N} \sum_{i=1}^N
Z_i' \Biggr| + \E \Biggl| \frac{1}{N} \sum
_{i=1}^N Z_i''
- \E\frac{1}{N} \sum_{i=1}^N
Z_i'' \Biggr| \\
&=&\!: E' +
E''.
\end{eqnarray*}
By Jensen's inequality, independence and the bound on $Z_i'$, we have
\[
\bigl(E'\bigr)^2 \le\Var \Biggl( \frac{1}{N}
\sum_{i=1}^N Z_i'
\Biggr) = \frac{1}{N^2} \sum_{i=1}^N
\Var\bigl(Z_i'\bigr) \le\frac{K^2}{N}.
\]
Moreover, by triangle and Jensen's inequalities,
\[
E'' \le2 \E \Biggl| \frac{1}{N} \sum
_{i=1}^N Z_i''
\Biggr| \le\frac{2}{N} \sum_{i=1}^N
\E\bigl|Z_i''\bigr|.
\]
The assumption on the tails of $Z_i$ implies that
$\P\{ |Z_i''| > t\} \le C/(t \vee K)^{1+\eta}$ for $t>0$, thus
\[
\E\bigl|Z_i''\bigr| = \int_0^\infty
\P\bigl\{ \bigl|Z_i''\bigr|>t \bigr\} \,dt \le
\frac{C}{K^\eta} + \frac{C}{\eta K^\eta} = C \biggl( 1 + \frac
{1}{\eta} \biggr)
K^{-\eta}.
\]
Hence
\[
E'' \le2C \biggl( 1 + \frac{1}{\eta} \biggr)
K^{-\eta}
\]
and
\[
E' + E'' \le\frac{K}{\sqrt{N}} + 2C
\biggl( 1 + \frac{1}{\eta} \biggr) K^{-\eta}.
\]
Choosing $K = (\varepsilon/2)\sqrt{N}$ and using the assumption on
$N$, one easily checks that
\[
E' + E'' \le\frac{\varepsilon}{2} +
\frac{\varepsilon}{2} \le \varepsilon
\]
as desired.
\end{pf*}

%apA #&#
\begin{appendix}\label{app}
\section*{\texorpdfstring{Appendix: Proof of Proposition \lowercase{\protect\ref{productdistributions}}}
{Appendix: Proof of Proposition 1.3}}

In this section we prove Proposition~\ref{productdistributions},
which states that product distributions satisfy the regularity assumption
in Theorem~\ref{main}.
Note that this result and its proof are not needed in the proof of
Theorem~\ref{main}.

Consider a random vector $X$ and an orthogonal projection $P$ in $\R^n$
as in Proposition~\ref{productdistributions}. Denoting by $(P_{ij})$ the
$n \times n$ matrix of the operator $P$, we express
\[
\|PX\|_2^2 = \langle X, PX\rangle= \sum
_{i,j=1}^n \xi_i \xi_j
P_{ij}.
\]
The contribution of the diagonal of $P$ to this sum is
\[
D:= \sum_{i=1}^n \xi_i^2
P_{ii}.
\]
Denote by $P_0$ the matrix $P$ with diagonal removed; then
%
%eA.1 #&#
\setcounter{equation}{0}
\begin{equation}
\label{removediagonal} \|PX\|_2^2 - D = \langle X,
P_0 X\rangle.
\end{equation}

We can estimate $\langle X, P_0 X\rangle$ using a standard
decoupling argument.
Let $X'$ denote an independent copy of $X$,
and let $\E_X$, $\E_{X'}$ denote the expectations with respect to $X$
and $X'$, respectively.
Since the matrix $P_0$ has zero diagonal, we have\footnote{Throughout
this proof, we write $a \lesssim b$
if $a \le C b$ for some constant $C$ which is independent of $n$.}
%
%eA.2 #&#
\begin{equation}
\label{decoupling} \E\bigl|\langle X, P_0 X\rangle\bigr|^p
\lesssim\E_{X'} \E_X \bigl|\bigl\langle X, P_0
X'\bigr\rangle\bigr|^p.
\end{equation}
This inequality can be obtained from general decoupling results; see
\cite{dlPG}, Theorem~3.1.1;
a simple and well-known proof of (\ref{decoupling}) is given in \cite
{Vdecoupling}.

Next, an application of a standard symmetrization argument and
Khintchine inequality
(or a direct application of Rosenthal's inequality~\cite{Ros},
see~\cite{FHJSZ}) yields for every $a \in\R^n$ that
\[
\E\bigl|\langle X, a\rangle\bigr|^p = \E \Biggl| \sum_{i=1}^n
a_i \xi_i \Biggr|^p \lesssim\|a
\|_2^p.
\]
Therefore, by conditioning on $X'$ we obtain from (\ref{decoupling}) that
%
%eA.3 #&#
\begin{equation}
\label{toP0X} \E\bigl|\langle X, P_0 X\rangle\bigr|^p \lesssim
\E_{X'} \bigl\|P_0 X'\bigr\|_2^p
= \E\|P_0 X\|_2^p.
\end{equation}

Since $P_0$ equals $P$ without the diagonal, the triangle inequality yields
\[
\|P_0 X\|_2 \le\|PX\|_2 + \Biggl( \sum
_{i=1}^n \xi_i^2
P_{ii}^2 \Biggr)^{1/2}.
\]
Since $0 < P_{ii} \le\|P\| \le1$, we can replace $P_{ii}^2$ by
$P_{ii}$, so
\[
\|P_0 X\|_2 \le\|PX\|_2 + D^{1/2}
\lesssim\bigl( \|PX\|_2^2 + D \bigr)^{1/2}.
\]
H\"older's inequality then implies that
%
%eA.4 #&#
\begin{equation}
\label{P0X} \E\|P_0 X\|_2^p \lesssim \bigl(
\E \bigl| \|PX\|_2^2 + D \bigr|^p \bigr)^{1/2}.
\end{equation}
Putting (\ref{removediagonal}), (\ref{toP0X}) and (\ref{P0X})
together, we
arrive at the inequality
\[
\E \bigl| \|PX\|_2^2 - D \bigr|^p \lesssim \bigl( \E \bigl|
\|PX\|_2^2 + D \bigr|^p \bigr)^{1/2}.
\]
Put in different words, the random variable $Z:= \|PX\|_2^2 - D$
satisfies the inequality
\[
\|Z\|_{L_p}^2 \lesssim\|Z + 2D\|_{L_p} \le\|Z
\|_{L_p} + 2\|D\|_{L_p}.
\]
Solving this quadratic inequality we obtain that
%
%eA.5 #&#
\begin{equation}
\label{Z} \|Z\|_{L_p} \lesssim1 + \|D\|_{L_p}^{1/2}.
\end{equation}

In order to bound $\|D\|_{L_p}$ we consider
\[
\|D-k\|_{L_p}^p = \E \Biggl| \sum_{i=1}^n
\xi_i^2 P_{ii} - k \Biggr|^p = \E \Biggl|
\sum_{i=1}^n \bigl(\xi_i^2
- 1\bigr) P_{ii} \Biggr|^p,
\]
where we\vspace*{1pt} used that $\sum_{i=1}^n P_{ii} = \tr(P) = k$.
Recall that by the assumptions we have $\E(\xi_i^2-1) = 0$ and
$\|\xi_i^2-1\|_{L_p} \le\|\xi_i^2\|_{L_p} + 1 = \|\xi_i\|_{L_{2p}}^2 + 1 \lesssim1$.
An application of Khintchine's inequality or Rosenthal's inequality (as
before) and the bound $P_{ii}^2 \le P_{ii}$ yield that
%
%eA.6 #&#
\begin{equation}
\label{D-kbound} \|D-k\|_{L_p}^p \lesssim \Biggl( \sum
_{i=1}^n P_{ii}^2
\Biggr)^{p/2} \le \Biggl( \sum_{i=1}^n
P_{ii} \Biggr)^{p/2} = \bigl(\tr(P)\bigr)^{p/2} =
k^{p/2}.
\end{equation}
It follows that
\[
\|D\|_{L_p} \le\|D-k\|_{L_p} + k \lesssim k^{1/2} +
k \lesssim k.
\]

Putting this into (\ref{Z}), we see that
%
%eA.7 #&#
\begin{equation}
\label{Zbound} \|Z\|_{L_p} \lesssim k^{1/2}.
\end{equation}
Finally, by definition of $Z$ and using the triangle inequality and
bounds (\ref{Zbound}), (\ref{D-kbound}), we conclude that
\[
\bigl\| \|PX\|_2^2 - k \bigr\|_{L_p} \le\|Z
\|_{L_p} + \|D-k\|_{L_p} \lesssim k^{1/2} +
k^{1/2} \lesssim k^{1/2}.
\]
Proposition~\ref{productdistributions} is proved.
\end{appendix}

\section*{Acknowledgments}
The authors are grateful to the referees whose comments improved the
presentation of the paper.\vadjust{\goodbreak}

%suskaldyti doi

% imsref loaded by lrinkeviciute, 2012-07-27 15:43:24
% imsref loaded by lrinkeviciute, 2012-07-27 15:52:51

\printaddresses

\end{document}